\documentclass[10pt,reqno]{article}
\usepackage[all]{xy}
\usepackage{amsmath}
\usepackage{amsthm}
\usepackage{amscd}
\usepackage{amssymb}

\input{epsf.tex}

\numberwithin{equation}{section}

\newcommand{\nb}[1]{#1\nobreakdash-}

\theoremstyle{definition}

\theoremstyle{plain}
\newtheorem*{theorem*}{Theorem}
\newtheorem{theorem}{Theorem}
\newtheorem{proposition}[theorem]{Proposition}
\newtheorem{lemma}[theorem]{Lemma}
\newtheorem{corollary}[theorem]{Corollary}

\newcounter{remarks}
{\paragraph*{Remarks}\smallskip
     \begin{list}{\arabic{remarks}. }{\usecounter{remarks}%
          \setlength{\leftmargin}{0in}%
          \setlength{\rightmargin}{0in}%
          \setlength{\labelsep}{0pt}%
          \setlength{\labelwidth}{0pt}%
          \setlength{\listparindent}{0pt}%
     }
}
{
\end{list}
}

\newcommand\inv\inverse

\DeclareMathOperator{\Ends}{Ends}
\DeclareMathOperator\QI{QI}
\DeclareMathOperator\QIhat{\widehat{QI}}
\DeclareMathOperator\interior{int}

\DeclareMathOperator\Stab{Stab}

\DeclareMathOperator\image{image}

\DeclareMathOperator\Isom{Isom}
\DeclareMathOperator\Homeo{Homeo}
\DeclareMathOperator\Ker{Ker}
\DeclareMathOperator\kernel{\Ker}

\DeclareMathOperator\Vertices{Verts}
\DeclareMathOperator\Edges{Edges}

\DeclareMathOperator\TIndex{TI}

\DeclareMathOperator\Star{Star}

\DeclareMathOperator\ad{ad}

\newcommand\R{{\mathbf R}}
\newcommand\reals{\R}
\newcommand\rational{{\mathbf Q}}
\newcommand\hyp{{\mathbf H}}
\newcommand\complex{{\mathbf C}}
\newcommand\Z{{\mathbf Z}}

\newcommand\E{{\mathcal E}}
\newcommand\C{{\mathcal C}}

\newcommand\Sum{\sum}
\newcommand\infinity{\infty}
\newcommand{\bdy}{\partial}
\newcommand\cobdy\delta
\newcommand{\from}{\colon}
\newcommand\composed{\circ}
\newcommand\suchthat{\bigm|}
\newcommand\inverse{{-1}}
\newcommand\union{\cup}
\newcommand\Union{\bigcup}
\newcommand\abs[1]{\left| #1 \right|}

\newcommand\subgroup{<}

\newcommand\homeo\approx
\newcommand\cbdy\delta

\newcommand\Id{{\text{Id}}}
\newcommand\intersect{\cap}

\newcommand\restrict{\bigm|}
\newcommand\semidirect{\rtimes}

\newcommand\cross{\times}

\newcommand\<\langle
\renewcommand\>\rangle
\newcommand\disjunion{\coprod}

\newcommand\V{{\mathcal V}}

\newcommand\G{{\mathcal G}}

\newcommand\Q{{\mathcal Q}}
\newcommand\tensor\otimes
\newcommand\VF{\mathcal{VF}}
\newcommand\floor[1]{\lfloor #1 \rfloor}
\newcommand\norm[1]{\left\| #1 \right\|}

\begin{document}

\setcounter{tocdepth}{2}

\title{Maximally symmetric trees}
\author{Lee Mosher,  Michah Sageev,  and Kevin Whyte}
\maketitle

\centerline{\itshape To John Stallings on his $64.8^{\text{th}}$ birthday}

\begin{abstract}
We characterize the ``best'' model geometries for the class of virtually
free groups, and we show that there is a countable infinity of distinct
``best'' model geometries in an appropriate sense--these are the maximally
symmetric trees. The first theorem gives several equivalent
conditions on a bounded valence, cocompact tree
$T$ without valence~1 vertices saying that $T$ is maximally symmetric. The
second theorem gives general constructions for maximally symmetric
trees, showing for instance that every virtually free group has a
maximally symmetric tree for a model geometry.
\end{abstract}

\section{Introduction} 

A \emph{model geometry} for a finitely generated group $G$ is a
proper metric space $X$ on which $G$ acts properly and coboundedly
by isometries; equivalently, there is a discrete, cocompact, finite kernel
representation $G \to \Isom X$. The group $G$ with its word metric is
quasi-isometric to each of its (coarse) geodesic model geometries (see
\S\ref{SectionCoarseGeodesic}).

Given a quasi-isometry class $\C$ of finitely generated groups, one can
ask:
\begin{enumerate}
\item Is there a common (coarse) geodesic model geometry $X$ for every
group in $\C$?
\item Is there a common locally compact group $\Gamma$, in which every
group of $\C$ has a discrete, cocompact, finite kernel
representation?
\end{enumerate}
A ``yes'' to (1) implies a ``yes'' to (2), with or without the word
``coarse''. We show in Corollary~\ref{CorollaryQuestions} that the coarse
version of (1) is equivalent to (2).

Many quasi-isometric rigidity theorems take the form of a positive answer
to these questions. For example, if $X$ is an irreducible nonpositively
curved symmetric space then $X$ is itself a model geometry for every group
quasi-isometric to $X$. For $X=\hyp^n$ see
\cite{Sullivan:ErgodicAtInfinity},
\cite{Tukia:quasiconformal}, \cite{CannonCooper}, and for $\complex\hyp^n$
see \cite{KoranyiReimann}, \cite{Chow:ComplexHyp}; in these cases, every
cobounded quasi-action on $X$ is quasiconjugate to an isometric action on
$X$. For $\rational\hyp^n$ and the Cayley hyperbolic plane see
\cite{Pansu:CC}, and for $X$ of higher rank see
\cite{KleinerLeeb:buildings} and also
\cite{EskinFarb}; in these cases, every quasi-isometry is a bounded
distance from an isometry.

By contrast, consider the class of groups $\VF$ which are virtually free
of finite rank $\ge 2$. $\VF$ is a single quasi-isometry class, as follows
from Stallings ends theorem \cite{Stallings:ends}, \cite{ScottWall} and
Dunwoody's accessibility theorem \cite{Dunwoody:Accessible}. $\VF$
coincides with the class of fundamental groups of finite graphs of finite
groups \cite{KPS:virtuallyfree}. The typical model geometries for $\VF$ are
bounded valence, bushy trees; any quasi-action on such a tree $T$ is
quasiconjugate to an isometric action on a \emph{possibly different} tree
$T'$ (\cite{MSW:QTOne}; see Theorem~\ref{TheoremTreeQuasiActions} below).
For the class $\VF$, questions (1), (2) have a negative answer
\cite{MSW:QTOne}: for primes $p\ne q\ge 3$ the groups $\Z/p * \Z/p$ and
$\Z/q * \Z/q$ have \emph{no} discrete, cocompact, virtually faithful
representations in the same locally compact group.

In lieu of a single model geometry for $\VF$, we describe ``best'' model
geometries $X$, namely those which are ``maximally symmetric''. Roughly
speaking this means that any continuous, proper, cocompact embedding of
$\Isom X$ into another locally compact group is an isomorphism. This is
asking for too much, though: one can always take the product of $X$
with a symmetric compact metric space; or equivariantly attach to
each point of some discrete orbit of $X$ a symmetric pointed compact
metric space. We should therefore avoid compact normal subgroups, and for
$\VF$ this is acheived by restricting to bounded valence, bushy trees with
no valence~1 vertices.  Our main results show how to recognize maximally
symmetric trees within this class, in quasi-isometric, topological, and
graph theoretical terms, and how to construct maximally symmetric trees
wherever needed; for example, each group in $\VF$ has some maximally
symmetric tree as a model geometry.

\paragraph{Statements of results} Quasi-isometries are reviewed in
\S\ref{SectionQI}. For any metric space $X$ the
\emph{quasi-isometry group} $\QI(X)$ is the group of self quasi-isometries
of $X$ modulo identification of quasi-isometries which have bounded
distance in the sup norm. Any quasi-isometry $f \from X \to Y$ induces an
isomorphism $\ad_f\from \QI(X) \to \QI(Y)$, and so any finitely
generated group with the word metric has the same quasi-isometry group as
any of its model geometries.

Let $X$ be a $\delta$-hyperbolic metric space for some $\delta \ge 0$. A
subgroup $H \subgroup \QI(X)$ is \emph{uniform} if its elements can be
represented by quasi-isometries of $X$ with uniform quasi-isometry
constants. Equivalently, $H$ can be represented by a quasi-action on $X$;
moreover, any two such quasi-actions differ in the sup norm by a bounded
amount. We say that $H$ is \emph{cobounded} if a representing quasi-action
is cobounded. 

A tree $T$ has \emph{bounded valence} if each vertex has valence $\le C$
for some constant $C$, and $T$ is \emph{bushy} if each vertex is a
uniformly bounded distance from some vertex $v$ such that at least three
components of $T-v$ are unbounded. Any two bounded valence, bushy trees
are quasi-isometric, and so their quasi-isometry groups are isomorphic. A
bounded valence, bushy tree $T$ is \emph{cocompact} if $\Isom T$ acts
cocompactly on $T$, or equivalently if the image of the natural
homomorphism $\Isom T \to \QI(T)$ is cobounded. A \emph{thorn} of
$T$ is a valence~1 vertex; if $T$ is thornless then $\Isom T$ has no
compact normal subgroups, and the homomorphism $\Isom T\to \QI(T)$ is
injective, among other nice properties.

\begin{theorem} [Characterizing maximally symmetric trees]
\label{TheoremMaximal}
For any bounded valence, bushy, cocompact, thornless tree $T$, the
following are equivalent:
\begin{itemize}
\item $\Isom T$ is a maximal uniform cobounded subgroup of $\QI(T)$.
\item For any bounded valence, bushy, thornless tree $T'$, any
continuous, proper, cocompact embedding $\Isom T \to \Isom T'$ is
an isomorphism.
\item 
For any locally compact group $\G$ without compact normal subgroups, any
continuous, proper, cocompact embedding $\Isom T \to \G$ is an
isomorphism. 
\end{itemize}
\end{theorem}

Such trees $T$ are called \emph{maximally symmetric}. 

The proof of Theorem~\ref{TheoremMaximal} is entirely abstract and
nonconstructive: it does not exhibit the existence of a single maximally
symmetric tree, let alone showing that they are model geometries for
virtually free groups. These facts follow from
Theorem~\ref{TheoremModelGeometries} whose proof is concrete and
constructive; as a byproduct, we obtain a finitistic method for
enumerating the isometry types of maximally symmetric trees. 

A bounded valence, bushy, cocompact tree $T$ is said to be
\emph{index~1 normalized} if, for every vertex $v$ and every incident edge
$e$, the stabilizers of $v$ and $e$ in $\Isom T$ satisfy
$[\Stab(v):\Stab(e)]\ge 2$. Note that index~1 normalized implies thornless.

\begin{theorem}[Existence of maximally symmetric trees]
\label{TheoremModelGeometries}
Fix a bounded valence, bushy tree $\tau$. 
\begin{enumerate}
\item Every uniform cobounded subgroup
of $\QI(\tau)$ is contained in a maximal uniform cobounded subgroup. 
\item For every maximal uniform cobounded subgroup $\G \subgroup
\QI(\tau)$ there exists a maximally symmetric tree $T$ which is index~1
normalized, and there exists a quasi-isometry $f \from T\to\tau$, such that
$\G = \ad_f(\Isom T)$. Moreover $T$ and $f$ are uniquely specified in the
following sense: if $\G = \ad_{f'}(\Isom T')$ for another
index~1 normalized $T'$ and quasi-isometry $f' \from T' \to \tau$, then
there exists an isometry $h\from T \to T'$ such that $\ad_f = \ad_{f'}
\composed \ad_h$. 
\item There is a natural one-to-one
correspondence between conjugacy classes of maximal uniform cobounded
subgroups of $\QI(\tau)$ and isometry classes of index~1 normalized,
maximally symmetric trees $T$. 
\item There is a countable infinity of such isometry classes; in fact
there is a countable infinity of both unimodular and nonunimodular
isometry classes. 
\end{enumerate}
\end{theorem}

Part (2) can be interpreted as saying that \emph{any} maximally symmetric
model geometry in the quasi-isometry class of $\tau$ is, in a certain
sense, equivalent to a maximally symmetric \emph{tree}. In part (4),
unimodularity of a tree $T$ means that the locally compact group $\Isom T$
is unimodular, i.e.\ each left invariant Haar measure on $\Isom T$ is also
right invariant. Unimodularity was shown by Bass and Kulkarni
\cite{BassKulkarni} to be equivalent to the existence of a discrete,
cocompact subgroup of $\Isom T$. 

\begin{corollary}
\label{CorollaryVirtuallyFree}
For every group $G \in \VF$ there exists a maximally symmetric tree $T$
which is a model geometry for $G$, and so $G$ is the fundamental group of a
finite graph of finite groups $\Gamma=T/G$ whose Bass-Serre tree $T$ is
maximally symmetric.
\end{corollary}

More information
about Corollary~\ref{CorollaryVirtuallyFree} can be extracted from the
proof of Theorem~\ref{TheoremModelGeometries}, namely an algorithm which
inputs any finite graph of finite groups, and outputs another one with the
same fundamental group whose Bass-Serre tree is maximally symmetric.

Our results should be compared and contrasted with results
concerning an irreducible, nonpositively curved symmetric space
$X$. For instance, recent results of Alex Furman
\cite{Furman:MostowMargulis} show that $X$ is ``maximally symmetric'' in a
manner very similar to that described in Theorem~\ref{TheoremMaximal}.
Contrasting with our Theorem~\ref{TheoremModelGeometries}, Furman's
results can be interpreted as saying that $X$ is the unique best
model geometry in its quasi-isometry class; also, the quasi-isometric
rigidity theorems quoted above show that $\QI(X)$ has a unique maximal
uniform cobounded subgroup up to conjugacy, namely $\Isom X$.

The proofs of Theorems~\ref{TheoremMaximal}
and~\ref{TheoremModelGeometries} use the rigidity theorem for
quasi-actions on trees from \cite{MSW:QTOne}, the theory of edge-indexed
graphs from \cite{Bass:covering} and \cite{BassKulkarni}, and tree
techniques reminiscent of \cite{BassLubotzky}. In particular, one the main
technical steps is Proposition~\ref{PropositionIsometry}, which gives
conditions on the quotient graphs of bounded valence, bushy, thornless
trees $T$, $T'$ that are sufficient to prove that any continuous, proper,
cocompact monomorphism $\Isom T \to \Isom T'$ is an isomorphism induced by
an isometry $T \to T'$. Proposition~\ref{PropositionIsometry} is related
to a result of Bass and Lubotzky (\cite{BassLubotzky}, Corollary 4.8(d)),
which obtains the same conclusion under somewhat stronger conditions.

The proof of Theorem~\ref{TheoremModelGeometries} gives a finitistic
characterization of maximally symmetric trees in terms of edge-indexed
graphs; see Corollary~\ref{CorollaryMaximal}. For example, the
bihomogeneous tree $T_{p,q}$ of alternating valences $p > q\ge 2$ is
maximally symmetric; the quotient edge-indexed graph $T_{p,q} /
\Isom T_{p,q}$ has two vertices and one edge, with one end of index $p$
and the other end of index $q$. In general, to each bounded valence, bushy
tree $T$ there corresponds a quotient edge-indexed graph $T/\Isom T$;
conversely to each edge-indexed graph $\Gamma$ there corresponds a
``universal covering tree'' $T$ and a deck transformation group $D(\Gamma)
\subgroup \Isom T$. With respect to this correspondence,
Corollary~\ref{CorollaryMaximal} describes a certain subclass of
edge-indexed graphs $\Gamma$ whose isomorphism classes are in one-to-one
correspondence with maximally symmetric trees. Also,
Corollary~\ref{CorollaryPumpingUpAlgorithm} gives a simple algorithm
which, given an edge indexed graph $\Gamma$, decides whether $\Gamma$
belongs to this subclass, and if not then the algorithm computes
another edge-indexed graph $\Gamma'$ which does belong, and for which
there is an nonsurjective embedding $D(\Gamma) \subgroup D(\Gamma')$ which
is continuous, proper, and cocompact. This algorithm can be used to prove
Corollary~\ref{CorollaryVirtuallyFree}, starting from a finite graph of
groups $\Gamma$ with fundamental group $G$. 

In the unimodular case, part (4) of Theorem~\ref{TheoremModelGeometries}
can be summarized by saying that there is a countable infinity of ``best''
geometries for the class~$\VF$. Part (4) is proved by simply giving some
examples, but we will improve that by showing in
Proposition~\ref{PropositionInfinite} that any finite graph $\Gamma$ in
which no edge is a loop and no two edges have the same endpoints has
infinitely many distinct edge-indexings, both unimodular and (when
$\Gamma$ is not a tree) nonunimodular, corresponding to a maximally
symmetric tree. The proof of Proposition~\ref{PropositionInfinite} will
give an effective, one-to-one enumeration of the isometry classes of
maximally symmetric trees.

While our focus in this introduction is mostly on the unimodular case,
Theorems~\ref{TheoremMaximal} and~\ref{TheoremModelGeometries} apply also
to nonunimodular trees. This may be applicable to graphs of groups having
bounded valence, bushy, Bass-Serre trees, in situations where these trees
can be nonunimodular, such as graphs of $\Z$'s, graphs of $\Z^n$'s, etc.

\subsubsection*{Acknowledgements}
\addcontentsline{toc}{subsection}{Acknowledgements}

The authors are supported in part by the National Science Foundation: the
first author by NSF grant DMS-9803396; the second author by NSF grant
DMS-989032; and the third author by an NSF Postdoctoral
Research Fellowship.

\section{Preliminaries}

\subsection{Metric spaces} 
\label{SectionMetricSpaces}
A metric space $X$ is \emph{proper} if closed balls are compact. This
implies that $X$ is complete, and that the isometry group $\Isom X$ is
locally compact and Hausdorff in the compact open topology. An action of a
group $G$ on $X$ will always mean an isometric action, that is, a
homomorphism $\phi \from G \to \Isom X$, usually written $g \mapsto\phi_g
\in \Isom X$. Properness of $X$ implies that an action is cocompact if and
only if it is cobounded, if and only if $\image(\phi)$ is a cocompact
subgroup of $\Isom X$ (a co-P action is one for which there is a P-subset
$K$ the union of whose translates equals the whole space). The action
$\phi$ is \emph{properly discontinuous} if for any two compact sets
$K,L \subset X$ the set $\{g \in G \suchthat \phi_g(K) \intersect L
\ne \emptyset\}$ is finite. When $X$ is proper, an action $\phi$ is
properly discontinuous if and only if $\image(\phi)$ is a discrete
subgroup of $\Isom X$ and $\kernel(\phi)$ is a finite subgroup of $G$.

\subsection{Quasi-isometries} 
\label{SectionQI}

A \emph{quasi-isometry} between two metric
spaces $X,Y$ is a map $f \from X \to Y$ such that for some constants $K \ge
1$, $C \ge 0$ we have
$$\frac{1}{K} d_X(x,y) - C \le d_Y(fx,fy) \le K d_X(x,y) + C, \quad x,y
\in X
$$
and for all $y \in Y$ there exists $x \in X$ such that $d_Y(fx,y) \le C$.
Every $K,C$ quasi-isometry $f \from X \to Y$ has a \emph{coarse inverse},
which is a $K,C'$ quasi-isometry $\bar f \from Y \to X$ such that
$d_{\sup}(\bar f \composed f, \Id_X) \le C'$ and $d_{\sup}(f \composed \bar
f,
\Id_Y) \le C'$, where the constant $C'$ depends only on $K,C$; the
notation $d_{\sup}$ denotes the sup metric on functions. 

The \emph{quasi-isometry group} of a metric space $X$, denoted $\QI(X)$,
is defined as follows. Let $\QIhat(X)$ denote the set of
quasi-isometries, equipped with the operation of composition. Define $f,g
\in \QIhat(X)$ to be \emph{coarsely equivalent} if $d_{\sup}(f,g) <
\infinity$, and let $[f]$ be the coarse equivalence class. Note that
composition is well-defined on coarse equivalence classes, thereby making
the set of coarse equivalence classes into a group
$\QI(X)$. The inverse of the coarse equivalence class of $f \in \QIhat(X)$
is the class of any coarse inverse for $f$.

Given a quasi-isometry $f \from X \to Y$ there is an induced isomorphism
$\ad_f \from \QI(X) \to \QI(Y)$ defined by $\ad_f[g] = [f \composed g
\composed \bar f]$ for any coarse inverse $\bar f$ of $f$.

A \emph{quasi-action} of a group $G$ on a metric space $X$ is a map $A
\from G \to \QIhat(X)$, denoted $g \mapsto A_g$, such that for some $K
\ge 1$, $C \ge 0$ we have: each map $A_g$ is a $K,C$ quasi-isometry;
$d_{\sup}(A_\Id,\Id) < C$; and for each $g,g' \in G$ we have $d_{\sup}(A_g
\composed A_{g'}, A_{gg'}) < C$. Postcomposing the quasi-action $G
\xrightarrow{A}\QIhat(X)$ with the quotient map $\QIhat(X)\to\QI(X)$ we
obtain the induced homomorphism $G\to\QI(X)$. The quasi-action $A$ is
\emph{cobounded} if there exists a bounded subset $D \subset X$ such that
for each $x \in X$ there is a $g \in G$ with $A_g(x) \in D$; also, $A$ is
\emph{proper} if for each $R>0$ there exists an integer $M > 0$ such that
for each $x,y\in X$ the cardinality of the set $\{g \in G \suchthat
d(A_g(x),y) \le R\}$ is at most $M$. Given quasi-actions $A,B$ of $G$ on
metric spaces $X,Y$ respectively, a \emph{quasiconjugacy} from $A$ to $B$
is a quasi-isometry $f\from X \to Y$ such that for some $C \ge 0$ we have
$d_{\sup}(f \composed A_g, B_g \composed f) \le C$, for all $g \in G$; it
follows that $\ad_f[A_g] = [B_g]$. Coboundedness and properness are
quasiconjugacy invariants of quasi-actions.

For quasi-isometries among hyperbolic metric spaces, boundary values
coarsely determine a quasi-isometry, in the following sense. For any
$\delta,K,C$ there exists $A$ such that if $f,g \from X \to Y$ are $K,C$
quasi-isometries between proper, geodesic, $\delta$-hyperbolic metric
spaces $X,Y$, and if the boundary extensions $\bdy f, \bdy g\from \bdy X
\to \bdy Y$ are identical, then $d_{\sup}(f,g)\le A$. It follows that when
$X$ is $\delta$-hyperbolic, the following two properties on a subgroup
$H\subgroup \QI(X)$ are equivalent: $H$ is \emph{uniform}, meaning that
each element of $H$ is represented by a $K,C$ quasi-isometry for some
fixed $K \ge 1$, $C \ge 0$; $H$ has an \emph{induced quasi-action}, namely
a quasi-action $s \from H \to\QIhat(X)$ such that the composition $H
\xrightarrow{s} \QIhat(X) \to \QI(X)$ equals the inclusion. When $H$ is
uniform, any two induced quasi-actions $s,s' \from H \to \QIhat(X)$ differ
by a bounded distance in the sup norm, that is,
$\sup\{d_{\sup}(sh,s'h)\suchthat h \in H\} <\infinity$. A uniform subgroup
$H\subgroup \QI(X)$ is \emph{cobounded} if some (and hence any) induced
quasi-action of $H$ on $X$ is cobounded.

\subsection{Coarse geodesic metric spaces}
\label{SectionCoarseGeodesic}

Sections~\ref{SectionCoarseGeodesic} and \ref{SectionLCGroups} contain the
proof of Corollary~\ref{CorollaryQuestions}, that the coarse
version of question~(1) in the introduction is equivalent to question~(2).
Beyond the basic definitions, most of the material of these two
subsections will not be needed for the rest of the paper.

In a metric space $X$, a \emph{geodesic} joining $x$ to $y$ is a path
$\alpha \from [a,b] \to X$ such that $x=\alpha(a)$,
$y=\alpha(b)$, and $d(\alpha(s),\alpha(t))=\abs{s-t}$ for $s,t \in [a,b]$.
We say that $X$ is a \emph{geodesic metric space} if any two points are
joined by a geodesic.

A \emph{coarse path} joining $x$ to $y$ is just a sequence
$x=x_0,\ldots,x_n=y$ in $X$; the \emph{word length} equals $n$, and the
\emph{path length} equals $\sum_{i=1}^n d(x_{i-1},x_i)$. We say that
$x_0,\ldots,x_n$ is a \emph{$C$-coarse path} if $d(x_{i-1},x_{i}) \le C$
for $i=1,\ldots,n$. A \emph{$C$-coarse geodesic} is a $C$-coarse path whose
path length equals the distance between its endpoints. A metric space $X$
is a \emph{coarse geodesic metric space} if there exists $C \ge 0$ such
that any two points are joined by a $C$-coarse geodesic.

A proper metric space $X$ is geodesic if and only if $d(x,y)$ is the
infimum of the path lengths of all rectifiable paths joining $x$ to $y$.
The next lemma, applied to the collection $\V$ of closed balls of radius
$C$, shows similarly that a proper metric space $X$ is $C$-coarse geodesic
if and only if $d(x,y)$ is the infimum of the path lengths of $C$-coarse
geodesics joining $x$ and $y$; we need a more general version of this
fact for later purposes.

We generalize the notion of a $C$-coarse path as follows. Let $\V=\{V(x)
\suchthat x \in X\}$ where for each $x$ the set $V(x) \subset X$ is a
compact neighborhood of $x$, and the following symmetry condition holds:
$x \in V(y)$ if and only if $y \in V(x)$. A \emph{$\V$-coarse path} is a
coarse path $x_0,\ldots,x_n$ such that $x_{i} \in V(x_{i-1})$ for
$i=1,\ldots,n$. We say that $X$ is \emph{$\V$-coarsely connected} if any
two points $x,y \in X$ can be joined by a $\V$-coarse path. In this case
we define the \emph{$\V$-word metric} $\mu_\V(x,y)$ to be the shortest word
length of a $\V$-coarse path joining $x$ to $y$, and the \emph{$\V$-path
metric} $\rho_\V(x,y)$ to be the infimum of the path lengths of all
$\V$-coarse paths joining $x$ to $y$.

\break

\begin{lemma}
\label{LemmaCoarseGeodesic}
Let $X$ be a proper metric space, and suppose that $\V=\{V(x)\}$ is as
above, and that $\V$ satisfies the following:
\begin{itemize}
\item $X$ is $\V$-coarsely connected.
\item There exists $R>r>0$ such that for each $x \in X$,
$$\overline B(x,r) \subset V(x) \subset \overline B(x,R)
$$
\end{itemize}
Then: $\rho_\V$ is a coarse geodesic metric whose restriction to each ball of
radius $r$ agrees with the given metric on $X$; and the metrics $\mu_\V$
and $\rho_\V$ are quasi-isometric, that is, the identity map is a
quasi-isometry between $\mu_\V$ and $\rho_\V$. 
\end{lemma}

\begin{proof} The idea of the proof is that $\rho_\V$ is a ``maximal
metric'' in the sense of Gromov \cite{Gromov:Asymptotic}, subject to the
constraint that $\rho_\V$ agrees locally with the given metric on $X$.

Given $x,y \in X$ and $n \ge \mu_\V(x,y)$, define
$\norm{x,y}_n$ to be the infimum of the path lengths of all $\V$-coarse
paths joining $x$ to $y$ which have word length $\le n$. Note that the
sequence
$\norm{x,y}_n$ is nonincreasing and has limit $\rho_\V(x,y)$. 

Fix $n$ for the moment. Since $X$ is proper, the infimum defining
$\norm{x,y}_n$ is acheived by some $\V$-coarse path $x=x_0,\ldots,x_k=y$
with a minimal word length $k=k(n)$, $\mu_\V(x,y) \le k \le n$.  Note that
the path $x_0,\ldots,x_k$ cannot have a subpath $x_{i-1},x_i,x_{i+1}$ such
that each of $d(x_{i-1},x_i)$, $d(x_i,x_{i+1})$ is $\le r/2$ because then
$d(x_{i-1},x_{i+1}) \le r$ which would produce a $\V$-coarse path
$x=x_0,\ldots,x_{i-1},x_{i+1},\ldots,x_k$ of path length $\le
\norm{x,y}_n$ whose word length is smaller than $k(n)$, a contradiction.
It follows that at least $\floor{k(n)/2}$ of the distances
$d(x_{i-1},x_i)$ are $>r/2$ (where $\floor{\bullet}$ denotes the greatest
integer function). We therefore have 
$$\norm{x,y}_n > (k(n)-1)C/4
$$

Now the sequence $k(n)$ is evidently nondecreasing; moreover,
$\norm{x,y}_{n+1} < \norm{x,y}_n$ if and only if $k(n) < k(n+1)=n+1$. If
$k(n)$ is not bounded above it follows that $\norm{x,y}_n$ diverges to
$+\infinity$, a contradiction. Therefore $k(n)$ is eventually constant,
proving that
$\norm{x,y}_n$ is eventually constant and equal to $\rho_\V(x,y)$. This
shows that $\rho_\V(x,y)$ is a coarse geodesic metric. 

Now we compare $\rho_\V$ to $\mu_\V$. Obviously
$$\rho_\V(x,y) \le R \cdot \mu_\V(x,y)
$$
For the other direction, we have seen that $\rho_\V(x,y)$ is realized by
some $\V$-coarse path of least word length
$k$, and the argument shows that
$$\mu_\V(x,y) \le k < \frac{4}{C} \rho_\V(x,y) + 1
$$
\end{proof}

\subsection{Locally compact, compactly generated groups} 
\label{SectionLCGroups}

All locally compact groups are assumed to be Hausdorff. For example,
from the Ascoli-Arzela theorem it follows that the isometry group of a
proper metric space is locally compact Hausdorff, in the compact open
topology.

The next lemma says that locally compact, compactly generated groups, like
finitely generated groups, have a well-defined geometry up to
quasi-isometry. Moreover, just as finite index implies quasi-isometry
among finitely generated groups, ``compact index'' implies quasi-isometry
among compactly generated groups.

\begin{lemma}
\label{LemmaCompactlyGenerated}
Let $\G$ be a locally compact topological group, $G$ a closed, cocompact
subgroup. Then $\G$ is compactly generated if and only if $G$ is compactly
generated. Moreover, if this is so then the inclusion $G \to \G$ is a
quasi-isometry with respect to the compactly generated word metrics.
Finally, any two compactly generated word metrics on $\G$ are
quasi-isometric.
\end{lemma}

\begin{proof} If we substitute ``finitely generated'' for ``compactly
generated'', and ``finite index'' for ``cocompact'', then this is a
standard result, and the proof goes through unchanged, with one caveat. To
show that two finite generating sets $A,B$ determine quasi-isometric word
metrics one must prove that $A \subset B^n$ and $B \subset A^m$ for some
integers $n,m$. We must prove the same when $A,B$ are compact generating
sets. 

First we reduce to the case of compact generating sets containing a
neighborhood of the identity $e$. Supposing that $A$ is any compact
generating set, it follows that $\union_{i=1}^\infinity A^n = \G$, and so
by the Baire category theorem some $A^i$ contains an open ball $B$. Also,
some $A^j$ contains $e$. Therefore, $A^{i+j}$ contains a neighborhood of
$e$, and we can replace $A$ by $A^{i+j}$.

Letting $A,B$ be two compact generating sets each containing a
neighborhood of $e$, for each $x \in B$ there exists $i$ such that $A^i$
contains a neighborhood of $x$, and by compactness of $B$ it follows that
$B \subset A^m$ for some $m$; similarly $A \subset B^n$.
\end{proof}

\paragraph{Remark on the proof} Note that $\reals$ has a compact generating
set with empty interior. Namely, if $E$ is any compact set with empty
interior and positive measure, then the set $E-E=\{e_1 - e_2 \suchthat
e_1,e_2 \in E\}$ contains a neighborhood of $0$, by an application of the
Lebesgue density theorem, and so $E\union -E$ generates $\R$ and has empty
interior.
\bigskip

The disadvantage of Lemma~\ref{LemmaCompactlyGenerated} is that a
compactly generated word metric does \emph{not} determine the correct
topology on $\G$: indeed, if the generating set contains a neighborhood of
the origin then the word metric is discrete. We correct this, at the same
time obtaining a coarse geodesic metric, as follows:

\begin{lemma}
\label{LemmaQuestions}
Suppose that $\G$ is a locally compact, compactly generated group. Then
there exists a left invariant coarse geodesic metric $\rho$ on $\G$ such
that $\rho$ yields the given topology on $\G$ and
$\rho$ is quasi-isometric to any compactly generated word metric on $\G$.
\end{lemma}

\begin{proof}
By a result of Birkhoff \cite{Birkhoff:NoteOnTopGroups} and of Kakutani
\cite{Kakutani:metrization}, we know that there exists a left invariant
metric $D$ on $\G$ yielding the topology on $\G$. Let $V$ be a compact
generating set for $\G$; by enlarging $V$ we may assume that $V$ contains
the $D$-ball of some radius $r>0$ about $e$, and that $V=V^\inv$. Let $\V =
\{g \cdot V \suchthat g \in G\}$. Applying Lemma~\ref{LemmaCoarseGeodesic},
the $\V$-coarse geodesic metric $\rho_\V$ agrees with $D$ on each $D$-ball
of radius $<r$, and $\rho_\V$ is quasi-isometric to the compactly
generated word metric $\mu_\V$. Moreover, $\rho_\V$ is clearly left
invariant.
\end{proof}

It follows that any group $G$ with a discrete, cocompact, finite kernel
representation to $\G$ preserves $\rho$ under the left action of $G$ on
$\G$, and so regarding $\G$ as a coarse geodesic, proper metric space $X$,
we obtain a properly discontinuous, cocompact action of $G$ on
$X$. This proves:

\begin{corollary} 
\label{CorollaryQuestions}
Given a collection of groups $\C$, the following are
equivalent:
\begin{itemize}
\item[(1)] There exists a coarse geodesic metric space $X$ on which each
group in $\C$ acts properly discontinuously and cocompactly.
\item[(2)] There exists a locally compact group $\G$ in which each group
of $\C$ has a discrete, cocompact, finite kernel representation.
\end{itemize}
\qed\end{corollary}

\subsection{Graphs}
\label{SectionGraphs} 

In this paper, \emph{all graphs and trees are locally finite}. Thus, a
\emph{graph} $\Gamma$ is a connected, locally finite
\nb{1}complex, and a \emph{tree} is a contractible graph. A \emph{vertex}
of $\Gamma$ means a \nb{0}cell; the set of vertices is denoted
$\Vertices(\Gamma)$. An \emph{edge} means a \nb{1}cell, that is, a
component of the complement of the vertices; the set of edges is denoted
$\Edges(\Gamma)$. For each edge
$e$ we choose a compact arc $\overline e$ and a characteristic map
$(\overline e,
\bdy\overline e) \to (\Gamma,\Vertices(\Gamma))$ taking
$\interior(\overline e)$ homemorphically to $e$. Each edge $e$ has two
\emph{ends} in the sense of Freudenthal, corresponding one-to-one with the
endpoints of $\overline e$, this correspondence being denoted $\eta
\leftrightarrow p_\eta$; the set of ends of $e$ is denoted $\Ends(e)$.
Each $\eta \in \Ends(e)$ is \emph{located} at a particular vertex of
$\Gamma$, namely the unique limit point in $\Gamma$ of the end $\eta$,
identified with the image of $p_\eta$ under the characteristic map.
Denote $\Ends(\Gamma) =\union_{e \in
\Edges(\Gamma)} \Ends(e)$. The set of $e \in \Ends(\Gamma)$ located at a
particular vertex $v \in \Vertices(\Gamma)$ is denoted
$\Ends(v)$; this corresponds to the ``link'' of $v$. We denote
$\Ends(e,v) = \Ends(e) \intersect \Ends(v)$, the set of ends of $e$
located at $v$, a set of cardinality zero, one, or two. An edge $e$ is
called a \emph{loop} if there exists $v \in \Vertices(\Gamma)$ such that
$\Ends(e,v) = \Ends(e)$. Note that we do not adopt a preferred orientation
for an edge, the distinction between the two orientations being encoded in
the two ends.

We impose on each graph $\Gamma$ a geodesic metric in which edge has
length~1. The \emph{isometry group} $\Isom \Gamma$ is defined to be the
group of cellular isometries of $\Gamma$; this coincides with the usual
isometry group except in the single case when $\Gamma$ is isometric to the
real line. Since $\Gamma$ is locally finite, it is proper, and so the group
$\Isom \Gamma$, with the compact-open topology, is locally compact
Hausdorff. Moreover, if $\Isom \Gamma$ acts cocompactly on
$\Gamma$ then $\Isom \Gamma$ is compactly generated: letting $\Delta$ be
any finite subgraph of $\Gamma$ whose translates under $\Isom\Gamma$ cover
$\Gamma$, the set $K_\Delta = \{g \in \Isom \Gamma \suchthat
g(\Delta) \intersect \Delta \ne \emptyset\}$ is a compact generating set.

\section{Characterizing maximally symmetric trees}
\label{SectionCharacterizing}

Given a bounded valence, bushy, thornless, cocompact tree $T$, to prove
Theorem~\ref{TheoremMaximal} we must prove the equivalence of the
following properties, which we may then take as the definition of
\emph{maximally symmetric}:
\begin{enumerate}
\item[(1)] 
For any locally compact group $\G$ with no compact normal subgroups, any
continuous, proper, cocompact monomorphism $\Isom T \to \G$ is an
isomorphism. 
\item[(2)] 
For any bounded valence, bushy, thornless tree $T'$, any
continuous, proper, cocompact monomorphism $\Isom T \to \Isom T'$ is
an isomorphism.
\item[(3)] 
$\Isom T$ is a maximal uniform cobounded subgroup of $\QI(T)$.
\end{enumerate}

\begin{proof}[Proof that (1) implies (2)] Obvious. \end{proof}

\begin{proof}[Proof that (2) implies (3)] 
Suppose that $\Isom T < A$ for some uniform cobounded subgroup $A$ of
$\QI(T)$. Choose an induced cobounded quasi-action $s \from A \to
\QIhat(T)$; as remarked at the end of \S\ref{SectionQI}, $s$ is unique up
to bounded distance in the sup norm. Now we apply the main result of
\cite{MSW:QTOne}:

\begin{theorem}[Rigidity of quasi-actions on trees]
\label{TheoremTreeQuasiActions}
If $T$ is a bounded valence, bushy tree and $s \from G \to \QIhat(T)$ is a
quasi-action of a group $G$ on $T$, then there exists an action $s' \from G
\to \Isom T'$ of $G$ on a bounded valence, bushy tree $T'$, and there
exists a quasiconjugacy $f \from T \to T'$ from $s$ to $s'$.
\qed\end{theorem}

We obtain a quasiconjugacy $f \from T\to T'$ from the quasi-action $s
\from A \to \QIhat(T)$ to an injective cobounded action $s' \from A \to
\Isom T'$. Restricting to $\Isom T$ gives an injective action $s' \from
\Isom T \to \Isom T'$ which is quasiconjugate via $f$ to the canonical
action of $\Isom T$ on $T$. Since $\Isom T$ is cobounded on $T$ it
follows that $s'(\Isom T) \subgroup\Isom T'$ is cobounded, that is,
cocompact, on $T'$.

Now we need a lemma from \cite{MSW:QTOne}:

\begin{lemma} 
\label{LemmaProper}
Given a bounded valence, bushy tree $T$, a sequence $(g_{i})$ converges in
$\Isom T$ if and only if $(g_i)$ satisfies the following property: 
\begin{description}
\item[Coarse convergence] There is a number $D$ so that
for any $v$ there is an $n$ so that the set $\{g_{i}(v) \suchthat i \geq 
n\}$ has diameter at most $D$.
\end{description}
\qed\end{lemma}

A convergent sequence $g_i \in \Isom T$ clearly satisfies coarse
convergence. Since coarse convergence is clearly invariant under
quasiconjugacy, the image sequence $s'(g_i) \in\Isom T'$ also satisfies
coarse convergence. Applying Lemma~\ref{LemmaProper} it follows that
$s'(g_i)$ converges in $\Isom T'$, proving that $s' \from
\Isom T \to \Isom T'$ is continuous. Also, $s'$ is proper, for suppose $C
\subset \Isom T'$ is compact. Choose a sequence $g_i \in s'{}^\inv(C)$.
Passing to a subsequence, $s'(g_i)$ converges to some $h \in C$. It follows
that $s'(g_i)$ satisfies coarse convergence in $T'$, and again by
quasiconjugacy invariance it follows that $g_i$ satisfies coarse
convergence in $T$. Applying Lemma~\ref{LemmaProper} it follows that $g_i$
converges in $\Isom T$ to some $g$. By continuity of $s'$ we have
$s'(g)=h$ and so $g \in s'{}^\inv(C)$, proving that $s'{}^\inv(C)$ is
compact and so $s'$ is proper.

Having proved that $s'$ is continuous, proper, and cocompact, applying (2)
it follows that $s'$ is surjective, which implies that $\Isom T=A$,
proving that $\Isom T$ is maximal in $\QI(T)$.
\end{proof}

\begin{proof}[Proof that (3) implies (1)] Assuming (3) is true suppose
that we have an embedding $\iota \from \Isom T \to \G$ as in (1). 

Let $\Delta$ be a compact fundamental domain for $T$ and consider the
compact generating set $K_\Delta = \{f \in \Isom T \suchthat f(\Delta)
\intersect \Delta \ne \emptyset\}$ for $\Isom T$. The left-invariant word
metric on $\Isom T$ determined by the generating set $K_\Delta$ is
quasi-isometric to the tree $T$. Specifically, the map $F\from
T\to\Isom T$, taking a vertex $w \in T$ to any isometry $F_w
\in\Isom T$ such that $w \in F_w(K_\Delta)$, is a quasi-isometry from $T$
to $\Isom T$. 

We have a quasi-isometry $F \from T \to \Isom T$, and applying
Lemma~\ref{LemmaCompactlyGenerated} the injection $\iota \from \Isom T
\to \G$ is a quasi-isometry. The left action of $\G$ on itself is
clearly a cobounded quasi-action, and quasiconjugating via $\iota
\composed F
\from T \to \G$ we
obtain a cobounded quasi-action of $\G$ on $T$. Applying
Theorem~\ref{TheoremTreeQuasiActions} produces a quasiconjugacy
$\Phi\from T \to T'$ from the $\G$ quasi-action on $T$ to a cobounded
action $A \from \G \to \Isom T'$ for some bounded valence, bushy
tree~$T'$. 

Repeating the argument above using Lemma~\ref{LemmaProper}, the
homomorphism $\G \xrightarrow{A} \Isom T'$ is continuous, proper,
and cocompact. Properness implies that the kernel is compact, but the group
$\G$ having no compact normal subgroup, it follows that $\G
\xrightarrow{A}
\Isom T'$ is an embedding. Letting $\bar \Phi \from T' \to T$ be a coarse
inverse of
$\Phi$, this shows that $\Isom T \subgroup
\ad_{\bar\Phi}(\Isom T')$, and the latter is clearly a uniform, cobounded
subgroup of $\QI(T)$. Applying (3) it follows that $\Isom T = \ad_{\bar
\Phi}(\Isom T')$, which implies that the composition of injections
$\Isom T\xrightarrow{\iota} \G \xrightarrow{A}
\Isom T'$ is an isomorphism, and so $\Isom T \xrightarrow{\iota} \G$ is
surjective.
\end{proof}

\section{Edge indexed graphs}
\label{SectionEdgeIndexedGraphs}

In this section we show how edge-indexed graphs can be used to encode
bounded valence trees. The material on graphs of groups and edge-indexed
graphs is taken for the most part from \cite{Bass:covering} and
\cite{BassKulkarni}. 

\subsection{Graphs of groups} For detailed references
see \cite{Serre:trees}, \cite{Bass:covering}, and \cite{ScottWall} for the
more topological viewpoint. We adopt a different notation for graphs than
these references. 

A \emph{graph of groups} is a graph $\Gamma$ together with a \emph{vertex
group} $\Gamma_v$ for each $v \in \Vertices(\Gamma)$, an \emph{edge group}
$\Gamma_e$ for each $e \in \Edges(\Gamma)$, and an \emph{edge-to-vertex
injection} $\gamma_\eta \from \Gamma_e \to \Gamma_v$ for each $\eta \in
\Ends(e,v)$. The fundamental group of $\Gamma$ is denoted $\pi_1\Gamma$,
and it acts on the Bass-Serre tree $T$. The definitions of $\pi_1\Gamma$,
of $T$, and of the action may be given topologically as in
\cite{ScottWall} or directly in terms of algebra as in
\cite{Serre:trees} or \cite{Bass:covering}, the link between the two
approaches being Van Kampen's theorem. Here is a brief account of the
topological definitions. 

For each vertex $v$ and edge $e$ choose a pointed, connected CW-complex
$X_v, X_e$ and an identification of the fundamental group $\pi_1 X_v$,
$\pi_1 X_e$ with the respective vertex or edge group $\Gamma_v, \Gamma_e$;
and for each end $\eta \in \Ends(e,v)$ choose a pointed cellular map
$\xi_\eta \from X_e\to X_v$ inducing the injection $\gamma_\eta$.
Construct a graph of spaces $X$ by gluing up the disjoint union of the
$X_v$'s and the products $X_e \cross \overline e$, where for each end
$\eta \in\Ends(e,v)$ we glue $X_e \cross p_\eta$ to $X_v$ via the gluing
$(x,p_\eta) \sim \xi_\eta(x)$ for each $x \in X_e$. For each vertex
$v$ of $\Gamma$ we define $\pi_1(\Gamma,v)$ to be $\pi_1(X,v_0)$ where $v_0
\in X_v$ is the base point. There is a natural quotient map $q \from X \to
\Gamma$, which induces a decomposition of $X$ into the point inverse
images $X_t=q^\inv(t)$, $t\in \Gamma$. Let $\tilde X$ be the universal
covering space of $X$. The components of lifts of decomposition elements
of $X$ defines a decomposition of $\tilde X$, and the corresponding
decomposition space of
$\tilde X$ is the tree
$T$. Choosing a base point $\tilde v \in \tilde X$ lying over $v \in X$
determines an identification of $\pi_1(\Gamma,v)$ with the deck
transformation group of the covering map $\tilde X \to X$, and the action
of $\pi_1(\Gamma,v)$ respects the decomposition of $\tilde X$ and so
descends to the required action of $\pi_1\Gamma$ on $T$.

In \S\ref{SectionUniversalCovering} we will review briefly the construction
of the Bass-Serre tree given in~\cite{Bass:covering}.

\subsection{Edge-indexed graphs} An
\emph{edge-indexing} of a connected graph $\Gamma$ is a function $I \from
\Ends(\Gamma) \to
\Z_+ =
\{1,2,3,\ldots\}$. The pair $(\Gamma,I)$ is
called an \emph{edge-indexed graph}. Given $v \in \Vertices(\Gamma)$,
define the \emph{valence} of~$v$, as usual, to be the cardinality of
$\Ends(v)$, and define the \emph{total index} of $v$ to be
$$\TIndex(v) = \Sum_{\eta \in \Ends(v)} I(\eta)
$$
When the edge-indexing $I$ is understood we will sometimes drop it from
the notation and simply say that $\Gamma$ is an edge-indexed graph.

For example, there is a forgetful functor which associates, to each graph
of groups $\Gamma$ having finite index edge-to-vertex injections, an
edge-indexing $I$ such that if $\eta \in \Ends(e,v)$ then $I(\eta) =
[\Gamma_v : \gamma_\eta(\Gamma_e)]$. In this example, $\TIndex(v)$ equals
the valence of any vertex $\tilde v$ of the Bass-Serre tree of $\Gamma$
such that $\tilde v$ lies over~$v$. 

To view this example in a slightly
different way, let the group $G$ act on a tree $T$ with quotient graph
$\Gamma=T/G$; we may view the quotient map $T \to \Gamma$ as a morphism of
graphs, taking vertices to vertices, edges to edges, and ends to ends, as
long as we \emph{first} subdivide any edge of $T$ which is inverted by
$G$. Given $\eta \in \Ends(e,v) \subset\Ends(\Gamma)$, choose $\tilde \eta
\in\Ends(\tilde e, \tilde v) \subset \Ends(T)$ lying over $\eta$, and
define $I(\eta) = [\Stab_G(v):\Stab_G(e)]$; note that $I(\eta)$ is
independent of the choice of $\tilde\eta$.

\subsection{Covering maps} Given an edge-indexed graph $\Gamma$, to take
an \emph{elementary subdivision} of $\Gamma$ means to choose a subset of
the edges of $\Gamma$, add a new vertex to the interior of each chosen
edge, and assign index~1 to each end incident to a new vertex; each new
vertex has valence~2 and total index~$2$. Thus, under elementary
subdivision, an edge can be subdivided into at most two edges. A general
\emph{subdivision} of $\Gamma$ is the result of a finite sequence of
elementary subdivisions; now an edge can be subdivided into an arbitrary
finite number of edges.

Let $\Gamma_1, \Gamma_2$ be edge-indexed graphs. A continuous, surjective
map $\mu\from\Gamma_1\to\Gamma_2$ is called a \emph{covering map} if there
exists a subdivision $\Gamma'_1$ of $\Gamma_1$ such that the following
holds:
\begin{description}
\item[Cellularity] $\mu$ is a cellular map from $\Gamma'_1$ to $\Gamma_2$,
taking
$\Vertices(\Gamma'_1)$ to $\Vertices(\Gamma_2)$, and taking each edge $e$
of $\Gamma'_1$ homeomorphically to an edge $\mu(e)$ of $\Gamma_2$. There is
therefore an induced map $\mu \from \Ends(\Gamma'_1) \to \Ends(\Gamma_2)$.
\item[Subdivision normalization] Given an edge $e$ of $\Gamma_1$ and a
vertex $v$ of $\Gamma'_1$ in the interior of $e$, if $e',e''$ are the edges
of $\Gamma'_1$ incident to $v$ then $\mu(e')=\mu(e'')$.
\item[Even covering] Each end $\eta \in \Ends(\Gamma_2)$ is \emph{evenly
covered by
$\mu$}, which means: letting $w \in \Vertices(\Gamma_2)$ be the vertex to
which $\eta$ is attached, for each $v \in \Vertices(\Gamma'_1)$ such that
$v \in
\mu^\inv(w)$, we have
$$I(\eta) = \sum_{\eta' \in \Ends(v) \intersect \mu^\inv(\eta)} I(\eta')
$$
\end{description}
Here are a few more properties which follow immediately from the
definitions:
\begin{description}
\item[Total index preserved] For each vertex $v$ of $\Gamma'_1$, we
have $\TIndex(v) = \TIndex(\mu(w))$. 
\item[Folding of subdivision vertices] Let $v$ be a subdivision vertex of
$\Gamma'_1$, and so $\TIndex(v)=2$ and $\TIndex(\mu(v))=2$. By
subdivision normalization it follows that $\mu(v)$ has valence~1, with one
incident end of index~2.
\end{description}
This last property, which derives from subdivision normalization, is
needed to avoid unnecessary subdivision.  As we shall see below in
Lemma~\ref{LemmaSubdivision} it follows in generic cases that the
subdivision $\Gamma'_1$ is, in fact, just an elementary subdivision of
$\Gamma_1$. 

Here are some examples. 

If there is a pair of vertices $v, w \in \Vertices(\Gamma)$ and a pair
of edges $e \ne e'$ each of whose two ends are attached respectively to
$v,w$, then there is a covering map which identifies $e$ to $e'$
homeomorphically, and leaves the rest of $\Gamma$ unchanged. We'll refer to
this covering map as \emph{collapsing a bigon}:
$$
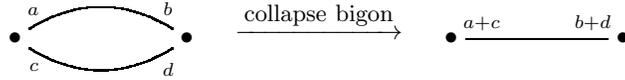
\hbox{\xymatrix{
\bullet \ar@/^1pc/@{-}^<<{a}^>>{b}[rr] \ar@/_1pc/@{-}_<<{c}_>>{d}[rr] &&
\bullet 
}}
\quad\xrightarrow{\text{\small collapse bigon}}\quad
\hbox{\xymatrix{
\bullet \ar@{-}[rr]^<<<{a+c}^>>>{b+d} && \bullet
}}
$$
This covering map is defined even when $v=w$.

If $e$ is a loop of $\Gamma$ then there is a covering map which first
subdivides $e$ and then collapses the resulting bigon, leaving the rest of
$\Gamma$ unchanged. This covering map
is called \emph{folding a loop}:
$$
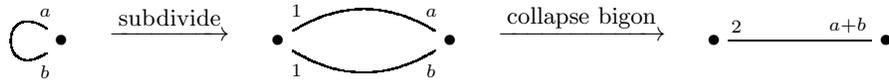
\hbox{\xymatrix{
\bullet \ar@(ul,dl)@{-}[]_<<{a}_>>{b}
}}
\quad\xrightarrow{\text{\small subdivide}}\quad
\hbox{\xymatrix{
\bullet \ar@/^1pc/@{-}^<<{1}^>>{a}[rr] \ar@/_1pc/@{-}_<<{1}_>>{b}[rr] &&
\bullet  
}}
\quad\xrightarrow{\text{\small collapse bigon}}\quad
\hbox{\xymatrix{
\bullet \ar@{-}[rr]^<<{2}^>>>>{a+b} && \bullet
}}
$$

If $\Gamma$ is a finite edge-indexed graph without loops or bigons, and if
any two vertices of $\Gamma$ have distinct total indices, it follows that
any covering map $\mu \from \Gamma \to \Gamma'$ is an isomorphism, because
$\mu$ must be one-to-one on vertices due to the fact that $\mu$ preserves
total indices, and the map on vertices determines the map on edges due to
the fact that $\Gamma$ has no loops or bigons.

Covering maps of edge-indexed graphs arise naturally from the covering
theory of graphs of groups \cite{Bass:covering}. Suppose that $\Gamma,
\Gamma'$ are graphs of groups with finite index edge-to-vertex injections,
and let $\Gamma,\Gamma'$ be equipped with their natural edge
indexings. Suppose that we have a covering map $\Phi \from\Gamma \to
\Gamma'$ in the graph of groups sense, as defined in \cite{Bass:covering}.
Then $\Phi$ is also a covering map in the edge-indexed graphs sense;
this follows from Proposition~2.7 of
\cite{Bass:covering}.

If $T$ is an arbitrary locally finite tree and $\G$ is an arbitrary
subgroup of $\Isom T$, then the quotient map $p \from T \to T/\G$ can be
regarded as a covering map with respect to a natural edge-indexing $I$ on
$T/\G$, where the edge-indexing on $T$ assigns index~1 to each end. To
describe $I$, first subdivide any edge of $T$ which is inverted by some
element of $\G$, so that $p$ is cellular. Given any vertex $w$ of $T/\G$
and end $\eta \in \Ends(w)$, choose $v \in p^\inv(w)$ and $\tilde\eta \in
p^\inv(\eta) \intersect \Ends(v)$, and define $I(\eta) = [\G_v : 
\G_{\tilde\eta}]$, where $\G_\bullet$ denotes a stabilizer subgroup of $\G$;
note that $I(\eta)$ is well-defined, independent of $v$ or $\tilde\eta$.
Note also that $I(\eta)$ is equal to the cardinality of the $\G_v$ orbit of
$\G_{\tilde\eta}$, which is useful in proving that $\eta$ is evenly covered
by $p$.

Here is a useful construction of covering maps:

\begin{lemma}
\label{LemmaCovering}
Let $T$ be a locally finite tree, and consider subgroups $\G \subgroup \G'
\subgroup \Isom T$ and the corresponding covering maps $T \xrightarrow{p'}
\Gamma' = T/\G'$ and $T\xrightarrow{p} \Gamma=T/\G$. The induced map
map $\Gamma' \xrightarrow{\mu} \Gamma$ is a covering map.
\end{lemma}

\begin{proof} We may assume, by subdividing $T$ if necessary, that $\G$
and $\G'$ act without edge inversions, and so the maps $p$, $p'$ are
cellular. Let $\G_\bullet,\G'_\bullet$ denote stabilizer subgroups of
$\G,\G'$ respectively. Consider a vertex $v$ of $T$, the image vertices
$w=p(v)$, $w'=p'(v)$, and $\eta \in \Ends(\Gamma,w)$; we must show that
$\eta$ is evenly covered by $\mu$. Let $E = p^\inv(\eta) \subset \Ends(v)$,
choose $\tilde\eta \in E$, and so the left hand side of the even covering
equation for $\eta$ is $[\G_v:\G_{\tilde\eta}]=\abs{E}$. Let
$\{\eta_1,\ldots,\eta_k\} = \mu^\inv(\eta) \intersect \Ends(w)$, and let
$E_i = p'{}^\inv(\eta_i) \intersect \Ends(v) = p'{}^\inv(\eta_i) \intersect
E$; choosing $\tilde\eta_i \in E_i$, the right hand side of the even
covering equation for $\eta$ equals the sum of
$[\G'_v:\G'_{\tilde\eta_i}]=\abs{E_i}$. But $E$ is the disjoint union of
$E_1,\ldots,E_k$.
\end{proof}

A finite edge-indexed graph $\Gamma$ is said to be an \emph{orbifold} if
each vertex has total index~2; it follows that topologically $\Gamma$ is
either a circle or an arc, and each vertex is either valence~2 with two
ends of index~1, or valence~1 with one end of index~2. For any covering
map
$\Gamma \to \Gamma'$ of edge-indexed graphs, $\Gamma$ is an orbifold if and
only if $\Gamma'$ is an orbifold, because of the fact that total index is
preserved. Covering maps between orbifolds can involve
complicated subdivisions. For example, if $\Gamma$ is a circle orbifold
and
$\Gamma'$ is an arc orbifold with one edge, first do \emph{any} subdivision
of $\Gamma$ resulting in an even number of edges, and then fold $\Gamma$
over $\Gamma'$ in zig-zag fashion. Similarly, if $\Gamma$ is an arc
orbifold and $\Gamma'$ is an arc orbifold with one edge, first do any
subdivision of
$\Gamma$ \emph{whatsoever}, and then fold the result in zig-zag fashion
over
$\Gamma'$. 

The following lemma demonstrates how subdivision normalization enforces
the simplest kind of subdivision for all but the most special covering
maps:

\begin{lemma}[Subdivision lemma]
\label{LemmaSubdivision}
If $p \from \Gamma_1 \to \Gamma_2$ is a covering map, and if $\Gamma_2$ is
not an arc orbifold with one edge, then the subdivision needed to
define $p$ is an elementary subdivision; in other words, each edge of
$\Gamma_1$ either maps homeomorphically to an edge of $\Gamma_2$ or is
folded around an edge of $\Gamma_2$.
\end{lemma}

\begin{proof}
Suppose that some edge $e$ of $\Gamma_1$ is subdivided by inserting at
least two distinct vertices in $\interior(e)$, and so $\Gamma'_1$ has an
edge $e'$ contained in the interior of $e$, with endpoints $a' \ne b' \in
\interior(e)$. Consider the edge $p(e')$ of $\Gamma_2$, whose ends are
located at vertices $p(a')$, $p(b')$. From the property ``folding of
subdivision vertices'' it follows that $p(a')$ and $p(b')$ both have
valence~1 and total index~2; this implies furthermore that $p(e')$ is the
unique edge of $\Gamma_2$.
\end{proof}

The next lemma satisfies one's natural intuition for covering
maps:

\begin{lemma} A composition of covering maps is a covering map.
\end{lemma}

\begin{proof} Consider a composition of covering maps  $\Gamma_1
\xrightarrow{\mu_1} \Gamma_2 \xrightarrow{\mu_2} \Gamma_3$. Let
$\Gamma'_1$ be the subdivision needed for $\mu_1$, and let $\Gamma'_2$ be
the subdivision needed for $\mu_2$. Pulling back the subdivision points of
$\Gamma'_2$ defines a further subdivision $\Gamma''_1$ of $\Gamma'_1$. The
map $\mu_2 \composed \mu_1$ from $\Gamma''_1$ to $\Gamma_3$ now satisfies
cellularity and subdivision normalization, and even covering is easily
checked.
\end{proof}

\subsection{The universal covering tree} 
\label{SectionUniversalCovering}

A \emph{universal covering map}
of an edge-indexed graph $\Gamma$ is a covering map $\pi \from T \to
\Gamma$ such that $T$ is a tree, regarded as an edge-indexed graph by
assigning index~1 to each end of each edge of $T$. Every edge-indexed
graph $\Gamma$ has a universal covering map. This is proved in Remark~1.18
of \cite{Bass:covering}; here is a construction.

Let $N_v$ denote the regular neighborhood of $v$ in $\Gamma$, defined as
the union of the \nb{1}cells of the barycentric subdivision of $\Gamma$
that touch $v$. The graph $N_v$ has the vertex $v$ and in addition one
valence~1 vertex denoted $m_\eta$ corresponding to each end $\eta \in
\Ends(v)$. Given $e \in \Edges(\Gamma)$ and $\eta \in \Ends(e)$ denote
$\eta^\inv \in \Ends(e)$ to be the end opposite from $\eta$, and note that
$m_\eta = m_{\eta^\inv}$. For each $v \in\Vertices(\Gamma)$ construct a
\emph{local universal cover} $p_v \from T_v \to N_v$, where $T_v$ is the
star on a set of cardinality $\TIndex(v)$, where $p_v$ is a
cellular map taking the star point to $v$, and where
$\abs{p_v^\inv(m_\eta)} = I(\eta)$ for each $\eta \in \Ends(v)$. Construct
$T$ and the covering map $p \from T \to \Gamma$ as the increasing union of
subtrees $T_0 \subset T_1 \subset T_2 \subset \ldots$ and maps $p_i \from
T_i \to\Gamma$, with $p_i \restrict T_j = p_j$ for $i>j$, as follows.
Choose a base vertex $v \in\Vertices(\Gamma)$, let $T_0$ be a disjoint
copy of $T_v$, and let $p_0$ be a disjoint copy of $p_v$. Assuming
$T_i, p_i$ have been constructed, consider an \emph{endpoint} $\tilde m$
of $T_i$, which means a valence~1 vertex of $T_i$ such that
$m=p_i(\tilde m)$ is not a vertex of $\Gamma$. Let $\tilde v$ be the vertex
of $T_i$ closest to $\tilde m$, let $v = p_i(\tilde v) \in
\Vertices(\Gamma)$, and note that $m=m_\eta$ for some $\eta\in \Ends(v)$.
The opposite end $\eta^\inv$ of $\eta$ is located at some vertex $w \in
\Vertices(\Gamma)$. Choose a disjoint copy of $T_w$, and choose a point
$m' \in p_w^\inv(m)$, a valence~1 vertex of $T_w$. Now glue the disjoint
copy of $T_w$ to $T_i$ by identifying $\tilde m$ to $m'$. Doing these
gluings disjointly for each valence~1 vertex $\tilde m$ of $T_i$ defines
the tree $T_{i+1}$, and extending $p_i$ by disjoint copies of the maps
$p_w$ defines the map $p_{i+1}$. This finishes the definition of the
universal covering tree $T$.

Note, following Remark~1.18 of \cite{Bass:covering}, that the Bass-Serre
tree of a graph of groups may be identified with the universal covering
tree of the underlying edge-indexed graph (this holds even when
edge-to-vertex injections are not of finite index, by stretching the
concept of an edge-indexing to accomodate arbitrary cardinal number
values for indices).

Starting from a finite valence tree $T$ and an action of a group
$G$ on $T$, take the graph of groups $T / G$, then pass to the associated
edge-indexed graph, and then take the universal covering tree; the result
is naturally isomorphic to the subdivision of $T$ obtained by elementarily
subdividing each edge which is inverted by some element of
$G$. Note in particular that if $T$ is not a line then the full metric
isometry group equals $\Isom T$, the group of cellular isometries, and so
in this case if $T'$ is the tree obtained by elementarily subdividing each
edge of $T$ that is inverted by some isometry of $T$ then
$\Isom T=\Isom T'$ and so $T'$ is the universal covering tree of $T' /
\Isom T'$.

\bigskip

We collect here without proof some simple facts, the first of which
justifies the terminology of a ``universal covering map'':

\begin{lemma} An edge-indexed graph $\Gamma$ with universal covering
$p\from T \to \Gamma$ satisfies the following properties: 
\begin{enumerate}
\item If $f\from\Gamma'\to\Gamma$ is a covering map and if $p' \from T'
\to \Gamma'$ is a universal covering map then there is an isomorphism
between $T$ and a subdivision of $T'$ so that $f \composed p' = p$.  
\item If $p'' \from \Gamma \to \Gamma''$ is a covering map then
the composition $T\xrightarrow{p}\Gamma \xrightarrow{p''} \Gamma''$ is a
universal covering map for $\Gamma''$.
\end{enumerate} 
\qed\end{lemma} 

It follows from (1) that a \emph{cellular} universal covering map $p \from
T\to\Gamma$ is uniquely determined by $\Gamma$ up to isomorphism: if
$p'\from T' \to \Gamma$ is another cellular universal covering map then
there is an isomorphism $\phi \from T \to T'$ such that $p' \composed \phi
= p$. The point here is that by definition a universal covering map $p
\from T\to\Gamma$ need not be cellular; there may be a nontrivial
subdivision in the definition of $p$.

\subsection{The geometric trichotomy} This is a simple trichotomy
satisfied by the universal covering tree $T$ of a finite, edge-indexed
graph $\Gamma$: 
\begin{description}
\item[$T$ is bounded;] or
\item[$T$ is line-like,] meaning that there is an embedded bi-infinite
line $L$ in $T$ and a constant $A \ge 0$ such that each point of $T$ is a
distance $\le A$ from $L$; or
\item[$T$ is bushy,] meaning that there is a constant $A \ge 0$ such that
each point of $T$ is a distance $\le A$ from a vertex $v$ with the
property that $T-v$ has at least three unbounded components.
\end{description}
When $T$ is line-like then $T$ has two ends; whereas when $T$ is bushy
then its space of ends is homeomorphic to a Cantor set. The proof of this
trichotomy is a standard exercise; see the comment before
Lemma~\ref{LemmaTrichotomy} below for an indication of a simple proof. See
\S 5 of \cite{BassKulkarni} for the statement and proof of the trichotomy
in the case when $\Gamma$ is unimodular; when $\Gamma$ is not
unimodular then it is easily seen that $T$ is bushy. As noted in
\cite{BassKulkarni}, in some sense this geometric trichotomy is yet
another manifestation of the spherical--euclidean--hyperbolic trichotomy. 

Note that the geometric trichotomy of the universal covering tree is
invariant under covering maps between finite edge-indexed graphs, because
if $p\from \Gamma_1 \to \Gamma_2$ is such a covering map then the universal
covering trees of $\Gamma_1, \Gamma_2$ are homeomorphic. 

The geometric trichotomy can be detected algorithmically from a finite
edge-indexed graph $\Gamma$ as follows. 

A \emph{thorn} of $\Gamma$ is a vertex $v$ with total
index~1; the valence must also equal~1. Equivalently, any vertex in the
universal covering tree $T$ lying over $v$ has valence~1. To \emph{trim} a
thorn means to remove it and the incident edge, producing a smaller
edge-indexed graph; the effect on the universal covering tree $T$ is to
remove the $D(\Gamma)$ orbit lying over $v$ and the incident edges (see
below for the definition of the deck transformation group $D(\Gamma)
\subgroup \Isom(T)$). Note that trimming does not affect the geometric
trichotomy of the universal covering tree. An edge-indexed graph
$\Gamma$ with no thorns is said to be
\emph{thornless}, and this happens if and only if the universal covering
tree $T$ has no valence~1 vertices, that is, if $T$ is also thornless.

Every finite edge-indexed graph $\Gamma$ can be trimmed inductively until
one reaches a thornless edge-indexed graph, which can be regarded as a
subgraph $\Gamma' \subset \Gamma$ called a \emph{thornless core} of
$\Gamma$. 

The following simple fact is left to the reader; its proof may be used to
provide a simple proof of the geometric trichotomy:

\begin{lemma} 
\label{LemmaTrichotomy}
Let $\Gamma$ be a finite edge-indexed graph with universal
covering $p \from T \to \Gamma$. Let $\Gamma'$ be a thornless core and let
$T' = p^\inv(\Gamma')$. Then $T'$ is a $D(\Gamma)$-invariant subtree of
$T$, and each point of $T$ is a uniformly bounded distance from some point
of $T'$. Moreover:
\begin{enumerate}
\item $T$ is bounded $\iff$ $T'$ is a single point $\iff$
$\Gamma'$ is a single point. 
\item $T$ is line-like $\iff$ $T'$ is a line $\iff$
$\Gamma'$ is an orbifold.
\item $T$ is bushy $\iff$ $T'$ has a vertex of valence~$\ge 3$ $\iff$
$\Gamma'$ is neither a single point nor an orbifold.
\end{enumerate}
In case (1) the thornless core may not be unique. In cases (2) and (3) the
thornless core $\Gamma'$ is unique.
\qed\end{lemma} 

Because of this lemma we may extend the terminology ``bushy'' to apply to
a finite edge-indexed graph $\Gamma$: we say that $\Gamma$ is bushy if and
only if its universal covering tree is bushy, which happens if and only if
$\Gamma$ has a unique thornless core $\Gamma'$, and $\Gamma'$ is neither a
point nor an orbifold. Many of the unpleasant properties we have begun to
encounter when the universal covering tree is a line may be simply avoided
by assuming bushiness.

It should be clear that a ``generic'' finite edge-indexed graph is bushy,
because generically the thornless core is neither a point nor an orbifold.

\subsection{Deck transformation groups} 
\label{SectionDeckTransformations}

Consider a finite, bushy edge-indexed graph $\Gamma$ with universal
covering $p\from T\to \Gamma$. We assume that $\Gamma$ has a geodesic
metric which lifts to a geodesic metric on $T$ so that each edge of $T$
has length~1. Recall that $\Isom T$ denotes the topological group of
cellular isometries of the tree $T$. The \emph{deck transformation group}
of $p \from T \to \Gamma$ is the closed subgroup of $\Isom T$ defined by
\begin{align*}
D(\Gamma) &= \{f \in \Isom T \suchthat p \composed f = p\} \\
          &= \{f \in \Homeo(T) \suchthat p \composed f = p\}
\end{align*}
The equation of sets on the right hand side is a consequence of the
Subdivision Lemma~\ref{LemmaSubdivision}, and it is the key place where
we need bushiness---the equation can fail when $\Gamma$ is an orbifold and
the covering map $p \from T \to \Gamma$ needs a nonelementary subdivision
on $T$. This equation implies that the topological quotient $T/D(\Gamma)$
is naturally identified with $\Gamma$, that is, there is a natural
homeomorphism
$T / D(\Gamma)\stackrel{\phi}{\approx} \Gamma$ so that the composition $T
\to T/D(\Gamma) \stackrel{\phi}{\approx}\Gamma$ equals~$p$. Moreover,
the homeomorphism $\phi$ induces an edge-indexed isomorphism between the
natural edge-indexed structure on $T/D(\Gamma)$ and $\Gamma$. 

To summarize: when $\Gamma$ is finite and bushy, the edge-indexed graphs
$T/D(\Gamma)$ and $\Gamma$ are naturally isomorphic. We will use this
identification without comment in what follows.

From finiteness of $\Gamma$ it follows that $D(\Gamma)$ is cocompact in
$\Isom T$.

\begin{lemma}
\label{LemmaTreeFacts}
Given a covering map of finite, bushy edge indexed graphs $\mu \from
\Gamma_1\to\Gamma_2$, with universal coverings $p \from T \to \Gamma_1$,
$\mu\composed p \from T \to \Gamma_2$, we have $D(\Gamma_1) \subset
D(\Gamma_2)$, with equality if and only if $\mu$ is an isomorphism of
edge-indexed graphs.

Given a finite, bushy edge-indexed graph $\Gamma$ with universal covering
$p\from T\to\Gamma$, and given a subgroup $\G \subgroup \Isom T$ with
$D(\Gamma)\subgroup \G$, the quotient map $T \to T/\G$ factors as $T
\xrightarrow{p} \Gamma \xrightarrow{p'} T/\G$ for some covering map $p'$.
\end{lemma} 

\begin{proof} The second part is simply a special case of
Lemma~\ref{LemmaCovering}. We leave the proof of the first part to the
reader, except to verify that if $\mu \from
\Gamma_1\to \Gamma_2$ is a nonisomorphic covering map then $D(\Gamma_1)
\neq D(\Gamma_2)$. To see why this is true, let $\Gamma'_1$ be the
subdivision of $\Gamma_1$, and note that there must two cells $c_1 \ne c_2$
of $\Gamma'_1$, either both edges or both vertices, such that
$\mu(c_1)=\mu(c_2)$. 
Choosing cells $\tilde c_1, \tilde c_2$ of $T$ lying over $c_1, c_2$
respectively, $D(\Gamma_2)-D(\Gamma_1)$ contains an isometry of $T$ taking
$\tilde c_1$ to $\tilde c_2$.
\end{proof} 

\begin{corollary}[Existence of a minimal subcover]
Let $\Gamma$ be a finite, bushy edge-indexed graph, and let $p \from T
\to \Gamma$ be the universal covering. There exists a covering map $\nu
\from \Gamma \to T / \Isom T$ which is a \emph{minimal subcover} for
$\Gamma$, meaning that for any covering map $\mu\from \Gamma
\to \Gamma'$, there exists a covering map $\mu' \from \Gamma' \to T /
\Isom T$ such that $\nu = \mu' \composed \mu$.

It follows easily that $D(\Gamma) = \Isom T$ if and only if $\Gamma$ is
its own minimal subcover.
\end{corollary}

\begin{proof} Consider the covering map $q \from T \to T / \Isom T$. Apply
Lemma~\ref{LemmaTreeFacts} to obtain a covering map factorization $T
\xrightarrow{p} \Gamma \xrightarrow{\nu} T/\Isom T$ of $q$. Consider any
covering map $\mu \from \Gamma \to \Gamma'$, and so the composition
$\mu\composed p\from T \to \Gamma'$ is a universal covering map. We have
$D(p) \subgroup D(\mu\composed p) \subgroup \Isom T$, which implies in
turn that $\nu\from\Gamma\to T /\Isom T$ factors as a product of covering
maps $\Gamma
\xrightarrow{\mu}\Gamma'\xrightarrow{\mu'} T / \Isom T$.
\end{proof}

To summarize this discussion, a finite, bushy edge-indexed graph $\Gamma$
can be considered as an encoding of a certain locally compact group,
namely the deck transformation group of the universal covering of
$\Gamma$. When $\Gamma$ is its own minimal subcover, then $\Gamma$ encodes
the entire isometry group of its universal cover. Thus, via the universal
covering map we obtain a bijection between the isometry types of
cocompact, bushy, bounded valence trees and the isomorphism types of
finite, bushy edge-indexed graphs which are their own minimal subcovers.

\subsection{Unimodularity} 

Given an edge-indexed graph $\Gamma$ there is a canonical cocycle $\xi \in
C^1(\Gamma;\rational_+)$, where $\rational_+$ denotes the group of positive
rational numbers under multiplication: for each oriented
edge $e$ of $\Gamma$ with positive end $\eta_+(e)$ and negative end
$\eta_-(e)$, define $\xi(e) = I(\eta_+(e)) / I(\eta_-(e))$.

\break 

\begin{theorem}[\cite{BassKulkarni}]
Let $T$ be a bounded valence, cocompact tree. The following are equivalent:
\begin{itemize}
\item $\Isom T$ has a discrete, cocompact subgroup.
\item $\Isom T$ is unimodular (that is, each left invariant Haar measure is
also right invariant).
\item The cohomology class of the canonical cocycle $\xi$ of $T/\Isom T$
is trivial.
\end{itemize}
\qed\end{theorem}

Because of this theorem, a finite edge-indexed graph $\Gamma$ is called
\emph{unimodular} if the canonical cocycle $\xi$ is cohomologically
trivial. For example if $\Gamma$ is a tree then it is unimodular.

\section{Pumping up the deck transformation group} 

Throughout this section and until further notice, \emph{all
edge-indexed graphs are finite and bushy}, allowing us to apply the
results of
\S\ref{SectionDeckTransformations}.

Let $\Gamma$ be a thornless edge-indexed graph with universal cover $T$.
We describe several methods for ``pumping up'' the group $D(\Gamma)$,
embedding it nonsurjectively in a larger deck transformation group. More
precisely, \emph{pumping up} $D(\Gamma)$ means constructing a continuous,
proper, cocompact, \emph{nonsurjective} monomorphism $\Psi\from D(\Gamma)
\to D(\Gamma')$ where $\Gamma'$ is another thornless edge-indexed graph. We
describe several explicit pumping up operations. One such operation has
already been discussed, namely ``passing to a proper subcover'', which
detects when the embedding $D(\Gamma) \subgroup \Isom(T)$ fails to be
surjective, and which allows one to construct another edge-indexed graph
whose deck transformation group \emph{is} all of $\Isom(T)$. Another
pumping up operation, called ``index~1 collapse'', may apply even when
$D(\Gamma) = \Isom(T)$, and it is particularly useful for demonstrating
that $T$ is not maximally symmetric. In addition we combine these two
operations into some composite pumping up operations. 

Given an edge-indexed graph $\Gamma$ which is its own minimal subcover, so
that $D(\Gamma) = \Isom(T)$, it follows that if $\Gamma$ can be pumped up
then $T$ fails to be maximally symmetric. Later on we will prove the
converse, giving the desired finitistic characterization of
maximally symmetric trees.

\subsection{Pumping up operations}

\paragraph{Proper subcovers} Passing to a proper subcover has already been
discussed: 

\begin{lemma}[Pumping up with proper subcovers] \quad For any proper
covering map $\mu \from \Gamma \to \Gamma'$, the inclusion map $D(\Gamma)
\subset D(\Gamma')$ (as subgroups of $\Isom T=\Isom T'$) is continuous,
proper, cocompact, and nonsurjective, and so $D(\Gamma)$ pumps up to
$D(\Gamma')$.
\qed\end{lemma}

\paragraph{Index~1 collapse} An \emph{index~1 edge} of $\Gamma$ is
an edge $e$ having an end $\eta$ of index~1; let $\zeta$ be the opposite
end of $e$. Suppose in addition that $e$ is not a loop 
(index~1 collapse on $e$ is \emph{not defined} when $e$ is a loop). The
ends $\eta,\zeta$ are located at distinct vertices $v \ne w$. Letting
$n = I(\zeta)$, we define the \emph{index~1--$n$ collapse
on~$e$}, also called an \emph{index~1 collapse} when the value of $n$ is
unimportant, to be the edge-indexed graph $\Gamma/e$ defined as follows.
The underlying graph of $\Gamma/e$ is obtained from $\Gamma$ by collapsing
the edge $e$ to a single vertex $z$. The quotient map $q\from\Gamma \to
\Gamma/e$ induces a bijection $\Ends(\Gamma)-\{\eta,\zeta\}
\leftrightarrow \Ends(\Gamma/e)$ denoted $\tilde\epsilon
\leftrightarrow\epsilon$. Define the index of each
$\epsilon \in \Ends(\Gamma/e)$ as follows:
$$
I(\epsilon) = \begin{cases}
I(\tilde\epsilon) & \text{if $\tilde\epsilon\not\in\Ends(v)$} \\
I(\tilde\epsilon) \cdot I(\zeta) = I(\tilde\epsilon) \cdot n & \text{if
$\tilde\epsilon\in\Ends(v)$}
\end{cases}
$$
For example:
$$
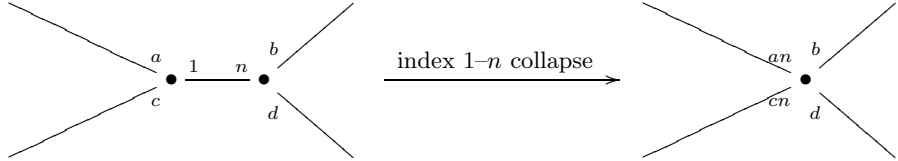
\hbox{\xymatrix{
\quad \ar@{-}[drr]^>>{a} & & & & \quad \ar@{-}[dl]_>>{b}  & & &
\quad \ar@{-}[drr]^>>{an\qquad} & & &
\quad \ar@{-}[dl]_>>{b}   \\
& & \bullet \ar@{-}[r]^<<{1}^>>{n} & \bullet & 
\quad \ar[rrr]^{\text{\small index 1--$n$ collapse}} & & &
& & \bullet    \\
\quad \ar@{-}[urr]_>>{c} & & & & \quad \ar@{-}[ul]^>>{d}  & & &
\quad \ar@{-}[urr]_>>{cn\qquad} & & &
\quad \ar@{-}[ul]^>>{d} 
}}
$$

Given an index~1--$n$ collapse $q \from \Gamma \to \Gamma/e$ and universal
coverings $p \from T\to\Gamma$ and $p'\from T'\to\Gamma/e$, the quotient
map $q$ lifts to a quotient map $\tilde q \from T \to T'$ that collapses
each connected component of $p^\inv(e)$ to a point, and $\tilde q$ is
equivariant with respect to a homomorphism $\Q\from D(\Gamma)\to
D(\Gamma/e)$. That is, for all $\phi \in D(\Gamma)$ we have a commutative
diagram
$$
\xymatrix{
T \ar[r]^{\tilde q} \ar[d]^{p} \ar@(ul,ur)^{\phi}
& T' \ar[d]^{p'} \ar@(ul,ur)^{\Q(\phi)}
\\
\Gamma \ar[r]^{q} & \Gamma/e
}
$$
We call $\Q$ the \emph{holonomy homomorphism} of the map $q$, and we note
that is a continuous, proper, cocompact monomorphism from $D(\Gamma)$ to
$D(\Gamma/e)$. 

Notice that $\tilde q$ is \emph{not} an isometry. On the other hand,
$\tilde q$ \emph{is} a quasi-isometry, because $\tilde q$ is one-to-one
except on the inverse image of each vertex $\tilde z$ of $T'$ lying over
$z=q(e) \in \Vertices(\Gamma/e)$, and the diameter of $\tilde
q^\inv(\tilde z)$ is at most~$2$.

\begin{lemma}[Pumping up with index 1--$n$ collapse] Let $\Gamma$ be a
thornless edge-indexed graph, and let $q \from \Gamma
\to \Gamma/e$ be an index~1--$n$ collapse with $p$, $p'$, $\tilde q$,
and $\Q$ as above.
If $n \ge 2$ then $\Q$ is not surjective and so $\Q$ pumps up $D(\Gamma)$
to $D(\Gamma/e)$. On the other hand, if $n=1$ then $\Q$ is an isomorphism.
\end{lemma}

\begin{proof} Let $v,w$ be the endpoints of $e$ incident to the ends
$\eta,\zeta$ of index $1,n$ respectively. Let $\tilde w \in T$ be a lift
of $w$. Let $e_1,\ldots,e_n \subset T$ be the incident lifts of
$e$, with ends $\zeta_1,\ldots,\zeta_n$ lifting $\zeta$ and incident to
$\tilde w$, and ends $\eta_1,\ldots,\eta_n$ lifting $\eta$ and incident to
vertices $\tilde v_1,\ldots,\tilde v_n$ lifting $v$, respectively. Since
$\Gamma$ is thornless and $I(\eta)=1$ it follows that there is an end
$\omega \in \Ends(v)-\{\eta\}$. Let $I(\omega)=m$. For each
$i=1,\ldots,n$, let $\Omega_i=\{\omega_{ij} \suchthat j=1,\ldots,m\}$, be
the lifts of $\omega$ located at $\tilde v_i$, and let $\Omega = \Omega_1
\union\cdots\union \Omega_n$. The subgroup $\Stab\Omega_i$ of $D(\Gamma)$
stabilizing $\Omega_i$ clearly acts as the symmetric group $S_m$ on the
set $\Omega_i$. The subgroup $\Stab\Omega$ acts on $\Omega$ preserving the
decomposition $\Omega = \Omega_1 \disjunion \cdots \disjunion \Omega_n$,
and so $\Stab\Omega$ acts on $\Omega$ as the semidirect product
$$\left(\underbrace{S_m \cross\cdots\cross S_m}_{\text{$n$ times}} \right)
\semidirect S_n
$$ 

Now consider the vertex $z = q(e) \in \Vertices(\Gamma/e)$ and its lift
$\tilde z = \tilde q(\tilde w) \in \Vertices(T')$. Note that the end
$q(\omega) \in
\Ends(z)$ has index
$mn$. Also, $\tilde q(\omega_{ij}) = \omega'_{ij} \in \Ends(\tilde z)$, and
letting $\Omega'_i = \{\omega'_{ij} \suchthat j=1,\ldots,m\}$ and $\Omega'
= \Omega'_1 \union\cdots\union \Omega'_n \subset \Ends(\tilde z)$, it
follows that $\Stab\Omega' \subset D(\Gamma/e)$ acts on $\Omega'$ as the
symmetric group $S_{mn}$. 

On the other hand, clearly $\Q(\Stab\Omega) \subset \Stab\Omega'$, and
the image of $\Q(\Stab\Omega)$ is the
subgroup of $\Stab\Omega'$ that preserves the decomposition $\Omega'
= \Omega'_1 \disjunion \cdots \disjunion \Omega'_n$. If $n \ge 2$ it
follows that
$\Q$ is not surjective, whereas when $n=1$ it follows that $\Q$ is
surjective.
\end{proof}

More generally, given a finite sequence of index~1 collapses
$\Gamma=\Gamma_1 \to\cdots\to \Gamma_n$, the resulting
composition has the effect of collapsing each component
of some forest $F\subset\Gamma$ of index~1 edges; we call this an
\emph{index~1 forest collapse} on $F$, denoted $\Gamma \xrightarrow{q_F}
\Gamma/F$, where $\Gamma/F\approx\Gamma_n$. Note that $q_F$ pumps up the
deck transformation group if and only if $F$ contains at least one edge of
index 1--$n$ with
$n \ge 2$. For every finite edge-indexed graph $\Gamma$ there exists a
\emph{maximal index~1 forest collapse} $\Gamma \xrightarrow{q_F} \Gamma'$,
resulting in an edge-indexed graph $\Gamma'$ which has no index~1 collapse;
equivalently, each index~1 edge of $\Gamma'$ is a loop. Note that it is
not possible to collapse an arbitrary forest of index~1 edges; for
example, if a vertex $v$ has two ends $\eta,\eta'$ each of index~1, lying
in edges $e,e'$ respectively, and if the opposite ends of $e,e'$ each have
index~$\ge 2$, then at most one of $e,e'$ is collapsed in any index~1
forest collapse. On the other hand, an arbitrary forest of index~1--1
edges can be collapsed. 

\paragraph{Collapse and subcover} Next we define a composite operation
called \emph{collapse and subcover}. Starting from an edge indexed graph
$\Gamma$, first do a maximal index~1 forest collapse, with the effect that
all remaining index~1 edges are loops; then pass to the minimal subcover,
which in particular has the effect of folding all loops. The resulting
graph $\Gamma'$ has no index~1 edges, and of course $\Gamma'$ is its own
minimal subcover. A collapse and subcover pumps up the deck transformation
group in either of the following two situations: an index~1--$n$ edge is
collapsed for some
$n \ge 2$; or the covering map is proper. A collapse and subcover fails to
pump up the deck transformation group, thereby inducing an isomorphism of
deck transformation groups, in the remaining case: all index~1 edges are
of type~1--1, they form a forest, and collapse of this forest produces an
edge-indexed graph which is its own minimal subcover.

For example, consider the edge-indexed graph
$$
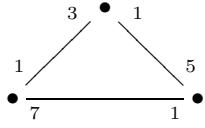
\xymatrix{
& \bullet \ar@{-}[dl]_>>{1}_<<{3} \ar@{-}[dr]^<<{1}^>>{5} \\
\bullet \ar@{-}[rr]_>>{1}_<<{7} && \bullet
}
$$
One can do an index~1 collapse on any one or two of the edges, but
not on all three. An index~1 collapse on
any two of these edges results in the edge-indexed graph
$$
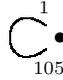
\xymatrix{
\bullet \ar@(ul,dl)@{-}[]_<<{1}_>>{105}
}
$$
and folding the loop results in the graph
$$
\xymatrix{
\bullet \ar@{-}[rr]^<<{2}^>>>{106} && \bullet
}
$$
which is its own minimal subcover.

\paragraph{Blowup and subcover} Because an index 1--1 forest collapse does
not change the deck transformation group, it is possible to start with an
edge-indexed graph which is its own minimal subcover, and then invert some
index 1--1 forest collapse, resulting in an edge-indexed graph which may
not be its own minimal subcover. We must therefore consider a composite
pumping up operation called ``blowup and subcover''.

Consider a thornless edge-indexed graph $\Gamma$. A \emph{blowup} of
$\Gamma$ consists of a thornless edge-indexed graph $\Gamma'$ which has a
forest $F$ of type 1--1 edges, and a collapse of this forest yielding
$\Gamma$, denoted $\Gamma\xleftarrow{q_F} \Gamma'$; formally $q_F$ has the
effect of identifying the collapsed graph $\Gamma'/F$ isomorphically
with~$\Gamma$. We also require that for each vertex $v$ of the blown up
graph $\Gamma'$, the valence and total index of $v$ are not both equal to
$2$; this, together with thornlessness of $\Gamma'$, has the important
implication that $\Gamma$ has only finitely many blowups up to
isomorphism. 

A \emph{blowup and subcover} of $\Gamma$, denoted $\Gamma\xleftarrow{q_F}
\Gamma' \xrightarrow{\mu} \Gamma''$, consists of a blowup $\Gamma
\xleftarrow{q_F}\Gamma'$ followed by a minimal subcover $\Gamma'
\xrightarrow{\mu} \Gamma''$; the resulting $\Gamma''$ is its own minimal
subcover. A blowup and \emph{proper} subcover pumps up the deck
transformation group, because the blowup $q_F$ induces an isomorphism
$D(\Gamma)\approx D(\Gamma')$, and the proper subcover $\mu$ induces a
nonsurjective monomorphism $D(\Gamma') \to D(\Gamma'')$. There are only
finitely many ways to blowup and subcover $\Gamma$, up to isomorphism.

For example, consider the edge-indexed graph 
$$
\xymatrix{
\bullet \ar@{-}[r]^<<{4}^>>{5} & \bullet \ar@{-}[r]^<<{3}^>>{6} & \bullet
}
$$
Note that this is its own minimal subcover, and so its deck transformation
group equals the isometry group of its universal cover $T$; however, $T$
turns out \emph{not} to be maximally symmetric, because there is a blowup
and proper subcover as follows. Blowing up the middle
vertex gives
$$
\xymatrix{
\bullet \ar@{-}[r]^<<{4}^>>{5} & 
   \bullet \ar@{-}[r]^<<{1}^>>{1} &
   \bullet \ar@{-}[r]^<<{3}^>>{6} & 
   \bullet
}
$$
Note that the total indices are 4, 6, 4, 6 from left to right. Folding
this graph up like a tri-fold wallet, 
$$
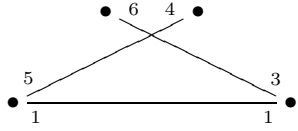
\xymatrix{
 & \bullet & \bullet \ar@{-}[dll]_<<{4}_>>{5} \\
\bullet \ar@{-}[rrr]_<<{1}_>>{1} &&& \bullet \ar@{-}[ull]_<<{3}_>>{6}
}
$$
we obtain a covering map to the
edge-indexed graph
$
\xymatrix{
\bullet \ar@{-}[r]^<<{6}^>>{4} & \bullet 
}
$.

A similar blowup and subcover may be carried out on any
edge-indexed graph of the form 
$$
\xymatrix{
\bullet \ar@{-}[rr]^<<<<{a+1}^>>{b} && \bullet \ar@{-}[rr]^<<{a}^>>>>{b+1}
&&
\bullet }
$$

\subsection{Pumping up algorithm} Is it possible to pump up $D(\Gamma)$
infinitely often? At first it may seem that one can indefinitely repeat
the blowup and proper subcover operation. However, we shall prove that
this is impossible (a fact which depends on bushiness), leading to an
algorithm which pumps up $\Gamma$ as much as possible with just a few
pumps.

Recall that we continue to assume all edge-indexed graphs are finite and
bushy.

\begin{lemma}
\label{LemmaAlgorithmStops}
Let $\Gamma$ be an edge-indexed graph without index~1 edges which is its
own minimal subcover. There is an algorithm which constructs a blowup and
subcover $\Gamma\xleftarrow{q_F}\Gamma'\xrightarrow{\mu}\Gamma''$ such that
$\Gamma''$ has no index~1 edge and no blowup and proper subcover. 
\end{lemma}

The algorithm of Lemma~\ref{LemmaAlgorithmStops} proceeds, in outline, as
follows. If $\Gamma$ has a blowup and proper subcover, then we show it has
one $\Gamma=\Gamma_0 \leftarrow \Delta_0 \rightarrow \Gamma_1$ such that
$\Gamma_1$ has no index~1 edges and is its own minimal subcover. Repeating
this we obtain a sequence of blowup and proper subcovers $\Gamma_{n-1}
\leftarrow \Delta_{n} \rightarrow \Gamma_{n}$ so that each $\Gamma_n$ has
no index~1 edges and is its own minimal subcover. This process stops at
some $\Gamma_N$ if and only if $\Gamma_N$ has no blowup and proper
subcover. To show that this eventually happens, we will prove that each
$\Gamma_n$ is connected to the original $\Gamma$ by a single blowup and
proper subcover operation $\Gamma\xleftarrow{F_n} \Gamma'_n \rightarrow
\Gamma_n$, and the number of edges in the collapsing tree $F_n$ is
increasing strictly monotonically with $n$. The crucial fact which makes
the algorithm stop is that at no stage does $\Gamma'_n$ ever have a vertex
of valence~2 and total index~2, which puts an upper bound on the number of
edges in $F_n$.

Lemma~\ref{LemmaAlgorithmStops} guarantees that the following algorithm
stops: 

\begin{corollary}[Pumping up algorithm] 
\label{CorollaryPumpingUpAlgorithm}
Given $\Gamma$ a finite, bushy, thornless
edge-indexed graph, the following algorithm pumps up $D(\Gamma)$ to
$D(\Gamma')$ where $\Gamma'$ is a finite, bushy
edge-indexed graph without index~1 edges and with no blowup and
proper subcover:
\begin{description}
\item[Step 1] Do a collapse and subcover $\Gamma \to \Gamma_1$,
and so $\Gamma_1$ has no index~1 edges and is its own minimal subcover.
\item[Step 2] If $\Gamma_1$ has no blowup and proper subcover, stop.
\item[Step 3] Otherwise, carry out the algorithm of
Lemma~\ref{LemmaAlgorithmStops} to find a blowup and proper subcover
$\Gamma_1 \leftarrow \Gamma_2 \rightarrow\Gamma_3$ such that $\Gamma_3$
has no index~1 edge and no blowup and proper subcover.
\end{description}
\qed\end{corollary}

\begin{proof}[Proof of Lemma~\ref{LemmaAlgorithmStops}] We break the
algorithm into two subroutines.

\begin{description}
\item[Subroutine 1] Given an edge-indexed graph $\Gamma$ which is its own
minimal subcover, and given any blowup and proper subcover
$\Gamma\xleftarrow{q_F}\Gamma_1\xrightarrow{\mu} \Gamma_2$, produce a
blowup and proper subcover $\Gamma \xleftarrow{q_{F'}} \Gamma' \to
\Gamma_4$ so that $\Gamma_4$ has no index~1 edges, and so that $F'
\ne\emptyset$.
\end{description}
Refer to the commutative diagram below. The fact that $F' \ne \emptyset$
follows because if not then $\Gamma \approx \Gamma'$ contradicting that
$\Gamma$ has no proper subcovers.

We may assume that $\Gamma_2$ does have at least one index~1 edge; note
that each of them is of type 1--1, because $\mu^\inv(e)$ for a type 1--$n$
edge
$e$ of $\Gamma_2$ is a union of type 1--$n$ edges of $\Gamma_1$, but each
index~1 edge of $\Gamma_1$ has type 1--1. Let $\Gamma_2
\xrightarrow{q_G}\Gamma_3$ collapse a maximal forest
$G$ of index 1--1 edges. Let $\Gamma_3\xrightarrow{\nu}\Gamma_4$ be a
minimal subcover, and so $\Gamma_4$ has no index~1 edge. Consider $\tilde
G=\mu^\inv(G)$, a subgraph of $F$ in $\Gamma_1$, and so $\tilde G$ is
itself a type 1--1 forest; in fact, $\mu$ induces an isomorphism between
each component of $\tilde G$ and a component of $G$. The type 1--1 forest
collapse $\Gamma_1\xrightarrow{q_F} \Gamma$ can be factored as a
composition of type 1--1 forest collapses $\Gamma_1\xrightarrow{q_{\tilde
G}} \Gamma' \xrightarrow{q_{F'}} \Gamma$ where $F'=F/\tilde G$; we remark
that no vertex of $F$ has valence and total index in $\Gamma_1$ both equal
to~2, and so the same is true of vertices of $F'$ in~$\Gamma'$. The
covering map $\mu \from \Gamma_1 \to \Gamma_2$ induces a covering map $\mu'
\from \Gamma' \to \Gamma_3$ so that $\mu' \composed q_{\tilde G} = q_G
\composed\mu$. Note that properness of $\mu$ implies properness of $\mu'$:
if some component of $G$ has more than one component in its preimage under
$\mu$ then this produces a point of $q_G(G)$ which has more than one
component in its preimage under $\mu'$; otherwise, some cell $c$ of
$\Gamma$ which is disjoint from $G$ has more than one preimage under
$\mu$, and it follows that $q_G(c)$ has more than one preimage under
$\mu'$. We thus obtain a blowup and proper subcover
$\Gamma\xleftarrow{q_{F'}}
\Gamma'\xrightarrow{\nu\composed \mu'}\Gamma_4$ so that $\Gamma_4$ has no
index~1 edges. 
$$
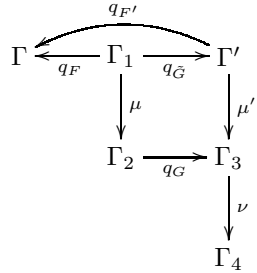
\xymatrix{
\Gamma & \Gamma_1 \ar[l]^{q_F} \ar[r]_{q_{\tilde G}} \ar[d]^{\mu} &
    \Gamma' \ar[d]^{\mu'} \ar@/_1pc/[ll]_{q_{F'}} \\
& \Gamma_2 \ar[r]_{q_G} & \Gamma_3 \ar[d]^{\nu} \\
& & \Gamma_4
}
$$
This completes the description of Subroutine~1.

\begin{description}
\item[Subroutine 2] Given successive blowup and subcover
operations $\Gamma\xleftarrow{q_F}\Gamma_1\xrightarrow{\mu}\Gamma_2$ and 
$\Gamma_2\xleftarrow{q_G} \Gamma_3 \xrightarrow{\nu} \Gamma_4$ so that
neither $\Gamma_2$ nor $\Gamma_4$ has any index~1 edge, produce another
blowup and subcover $\Gamma \xleftarrow{q_{F'}} \Gamma' \rightarrow
\Gamma_4$. Moreover, $\abs{F'} \ge \abs{F} + \abs{G}$.
\end{description}

First we consider the case where $G=\{e\}$, in which case we do a pushout
as follows:
$$
\xymatrix{
\Gamma & \Gamma_1 \ar[l]_{q_F} \ar[d]^{\mu} &
    \Gamma' \ar[d]^{\nu'} \ar[l]_{q'} \\
& \Gamma_2 & \Gamma_3 \ar[l]^{q_e} \ar[d]^{\nu} \\
& & \Gamma_4
}
$$
As a subset of $\Gamma_1 \cross \Gamma_3$, $\Gamma'$ is the set of ordered
pairs $(x,y)$ such that $\mu(x)=q_e(y)$, and the maps $q'$, $\nu'$ are
projections.  We denote the vertices and edges of $\Gamma'$ as
``tensor products''. Each vertex of
$\Gamma'$ is of the form $V \tensor W = (V,W)$ where $V\in
\Vertices(\Gamma_1)$, $W \in \Vertices(\Gamma_3)$, and
$\mu(V)=q_e(W)$. There are two types of edges in $\Gamma'$. First, letting
$Z=q_e(e) \in \Vertices(\Gamma_2)$, for each vertex $\tilde Z \in
\Vertices(\Gamma_1)$ such that $\mu(\tilde Z)=Z$ there is an edge $\tilde Z
\tensor e = \tilde Z \cross e\in \Edges(\Gamma')$. Second, for each edge $D
\in \Edges(\Gamma_1)$ and $E \ne e \in \Edges(\Gamma_3)$ such that
$\mu(D)=q_e(E)$, there is an edge $D \tensor E = \{(x,y) \in D
\cross E \suchthat \mu(x)=q_e(y)\}$. Edge indexing on $\Gamma'$ is
defined as follows: each edge of the form $\tilde Z \tensor e$ has both
ends of index~1; and the indexing of every other edge $D \tensor E$ is
obtained by pullback under the projection $D \tensor E\to
D \in \Edges(\Gamma_2)$.

We first check that $\nu' \from \Gamma' \to \Gamma_3$ is a covering map,
and the only thing to verify is that each end $\eta \in \Ends(\Gamma_3)$
is evenly covered. Noting that $\nu'{}^\inv(e)$ is the disjoint union of
the edges $\tilde Z \tensor e$, it follows that each end of $e$ is evenly
covered. If $\eta$ is not an end of $e$, let $\eta \in \Ends(W)$ for $W
\in \Vertices(\Gamma_3)$, and consider a vertex $V \tensor W \in
\nu'{}^\inv(W)$. Note that $I(\eta) = I(q_e(\eta))$. Note also that $q'$
induces an index preserving bijection between the set $\nu'{}^\inv(\eta)
\intersect \Ends(V \tensor W)$ and the set $\mu^\inv(q_e(\eta)) \intersect
\Ends(V)$. Since $q_e(\eta)$ is evenly covered, it follows that $\eta$ is
evenly covered.

Next we check that $q' \from \Gamma' \to \Gamma_1$ is a blowup of
$\Gamma_1$. Setting $F'$ to be the set of all edges of $\Gamma'$ of the
form $\tilde Z \tensor e$, note that these edges are pairwise disjoint and
so form an index 1--1 forest, and $q'$ is the map which collapses each
edge $\tilde Z \tensor e$ to the point $\tilde Z$. The only important issue
to resolve is whether $\Gamma'$ has a vertex of valence~2
and total index~2. Suppose there is such a vertex $V \tensor W$. Since
$\Gamma \xleftarrow{q_F} \Gamma_1$ is a blowup of $\Gamma$ it follows that
$\Gamma_1$ has no vertex of valence~2 and total
index~2, which implies that one of the two edges of $\Gamma'$
incident to $V\tensor W$ is an edge $\tilde Z \tensor e$ of $F'$, and so
$V=\tilde Z$ and $W$ is an endpoint of $e$. The two ends of $\Gamma'$
incident to $\tilde Z\tensor W$ are mapped
distinctly to $\Gamma_3$, and since we've already proved that $\nu'$ is a
covering map it follows that $W$ has valence~2 and total index~2 in
$\Gamma_3$. But $q_e \from\Gamma_3 \to\Gamma_2$ is a
blowup of $\Gamma_2$, which implies that $\Gamma_3$ has no vertices of
valence~2 and total index~2, a contradiction.

Finally, we must check that the composition $q_F \composed q' \from
\Gamma' \to \Gamma$ is a blowup of $\Gamma$. Clearly $F''=q'{}^\inv(F)
\union F'$ is an index~1--1 forest in $\Gamma'$, and the composition $q_F
\composed q'$ is just collapsing of $F''$. And we have already checked
above that $\Gamma'$ has no vertex of valence~2 and total index~2. Note
also that $\abs{F''} = \abs{q'{}^\inv(F)} + \abs{F'} =
\abs{F} + \abs{F'} \ge \abs{F} + \abs{G}$.

We therefore have the required blowup and subcover $\Gamma
\leftarrow \Gamma' \rightarrow \Gamma_4$ when $G$ is single edge. More
generally we can proceed inductively, doing successive pushouts, to obtain
the blowup and subcover for general $G$. This completes the description of
Subroutine~2.

\bigskip

To prove Lemma~\ref{LemmaAlgorithmStops}, consider an edge-indexed graph
$\Gamma$ without index~1 edges which is its own minimal subcover. If
$\Gamma$ has no blowup and proper subcover, the algorithm is finished.
Otherwise, choose a blowup and proper subcover $\Gamma\xleftarrow{q_{F_1}}
\Gamma'_1 \xrightarrow{\mu_1} \Gamma_1$, and immediately apply
Subroutine~1 to obtain one for which $\Gamma_1$ has no index~1 edges and
$F_1 \ne \emptyset$. Note that $\Gamma_1$ is its own minimal subcover. Now
proceed inductively as follows: assume we have a sequence of blowup and
proper subcover operations $\Gamma
\xleftarrow{q_{F_k}}\Gamma'_k \xrightarrow{\mu_k} \Gamma_k$ for which
$\Gamma_k$ has no index~1 edges and is its own minimal subcover, and such
that the cardinalities $\abs{F_k}$ are strictly increasing. If $\Gamma_k$
has no blowup and proper subcover then the algorithm is finished.
Otherwise, choose a blowup and proper subcover $\Gamma_k
\xleftarrow{q_{G_k}} \Delta_{k+1} \xrightarrow{\nu_k} \Gamma_{k+1}$, and
immediately apply Subroutine~1 to obtain one for which
$\Gamma_{k+1}$ has no index~1 edges and $G_k \ne \emptyset$. Now apply
Subroutine~2 to obtain a blowup and proper subcover $\Gamma
\xleftarrow{q_{F_{k+1}}}
\Gamma'_{k+1} \xrightarrow{\mu_{k+1}} \Gamma_{k+1}$, and note that
$\abs{F_{k+1}} \ge \abs{F_k} + \abs{G_k} > \abs{F_k}$. As the
cardinalities of the blowup trees $F_k$ are increasing, the blowup and
subcovers
$\Gamma\xleftarrow{q_{F_k}} \Gamma'_k \xrightarrow{\mu_k} \Gamma_{k}$ are
all distinct. But $\Gamma$ has only finitely many blowup and subcovers,
and so the algorithm must stop.
\end{proof}

\section{Enumerating maximally symmetric trees}
\label{SectionEnumerating}

In this section we prove Theorem~\ref{TheoremModelGeometries} and
Corollary~\ref{CorollaryVirtuallyFree}. For the
proof we will use without comment the fact that a bounded valence, bushy,
cocompact tree is index~1 normalized if and only if its quotient
edge-indexed graph has no index~1 edges. We continue to assume that all
edge-indexed graphs are finite and bushy. 

\subsection{Graphs with no blowup and proper subcover}

Besides applying Theorem~\ref{TheoremTreeQuasiActions}, the heart of the
proof is the following:

\begin{proposition}
\label{PropositionIsometry}
Let $\Gamma$ be an edge-indexed graph with no index~1 edge and no blowup
and proper subcover, and let $p \from T \to \Gamma$ be the universal
covering. Let $\Gamma'$ be an edge-indexed graph with no index~1 edge, and
let $p \from T' \to \Gamma'$ be the universal covering. If $\Psi\from\Isom
T \to \Isom T'$ is a continuous, proper, cocompact monomorphism, then
there exists an isometry $\psi \from T \to T'$ such that $\Psi =\ad_\psi$.
\end{proposition}

The techniques of proof are similar to those of Bass and Lubotzky
\cite{BassLubotzky}, who study situations under which a morphism between
two actions of a group $H$ on trees $T,T'$ is actually an isometry between
$T$ and $T'$. In the context of Proposition~\ref{PropositionIsometry}, if
one assumes in addition that $\Gamma$ has all ends of index~$\ge 3$, and
that the edge-indexed graph $T' / \image(\Psi)$ has all ends of index~$>
1$, then the conclusion follows from \cite{BassLubotzky}~Corollary 4.8(d).
Proposition~\ref{PropositionIsometry} makes no assumptions about how
$\Psi(\Isom T)$ acts on $T'$, other than the mild assumptions of
continuity, properness, and cocompactness of $\Psi$. Even more
significant, new techniques are needed in order to handle index~$2$ ends
of~$\Gamma$. 

\begin{proof} Let $p \from T \to \Gamma$, $p' \from T' \to \Gamma'$, $\Psi
\from \Isom T \to \Isom T'$ be as in the statement of the proposition. We
must produce an isometry $\psi \from T \to T'$ such that $\Psi=\ad_\psi$,
or in other words $\psi$ is \emph{$\Psi$-equivariant} meaning that for any
$f\in \Isom T$ we have $\psi \composed f = \Psi(f) \composed \psi$. Keeping
in mind the results of Bass and Lubotzky \cite{BassLubotzky} mentioned
above, our main difficulties are to understand index~2 ends of $\Gamma$.

Note that since $\Gamma$ and $\Gamma'$ are their own minimal
subcovers, we have $D(\Gamma) = \Isom T$, $D(\Gamma')=\Isom T'$.

We collect some facts about the tree $T$. Notation: the subgroups of
$\Isom T$ stabilizing a vertex $v$ and an edge $e$ are denoted $S_v$,
$S_e$ respectively. 

We may assume that $\Isom T$ (and similarly $\Isom T'$) acts without edge
inversions, and so for each vertex $v$ and each incident edge $e$ we have
$S_e \subgroup S_v$; if an edge is inverted by some isometry,
simply subdivide the edge at its midpoint.

For each vertex $v \in \Vertices(T)$ and each edge $e$ of $T$ with
endpoint $v$ and end $\eta$ incident to $v$, we have
$[S_v:S_e] = I(p(\eta))$. The group $S_v$ acts on the edges incident to
$v$, and it follows that $I(p(\eta))$ equals the cardinality of the
$S_v$-orbit of $e$.

An edge $e$ of $T$ with endpoints $v,w$ is called an \emph{index~2 edge}
if $S_e$ has index~2 in at least one of $S_v$ or $S_w$, or equivalently if
the image of $e$ in $\Gamma$ has an index~2 end. 

\begin{lemma}
\label{LemmaIndexTwo}
Given two edges $e,e'$ of $T$, the following are equivalent:
\begin{itemize}
\item[(a)] One of $S_e, S_{e'}$ is a subgroup of the other.
\item[(b)] $S_e = S_{e'}$.
\item[(c)] Letting $e=e_0 * e_1 * \cdots * e_k = e'$ be the unique
embedded edge path in $T$ from $e$ to $e'$, and letting $v_i = e_{i-1}
\intersect e_i$, for each $i=1,\ldots,k$ the set $\{e_{i-1},e_i\}$ forms a
single orbit (of cardinality~2) for the action of $S_{v_i}$ on
set $\Edges(v_i)$ of edges incident to $v_i$.
\end{itemize}
\end{lemma}

\begin{proof}
Obviously (b) implies (a). 

To show that (c) implies (b) it suffices to
observe that if $e \ne e'$ are both incident to a vertex $v$ and if
$\{e,e'\}$ forms an orbit of the action of $S_v$ then $S_e = S_{e'}$. 

To prove that (a) implies (c), suppose that $S_e \subset S_{e'}$ and let
$e=e_0*\cdots*e_k=e'$ be the edge path as in (c). It follows that $S_e
\subset S_{e_i}$ for $i=1,\ldots,k$. By induction on $k$ we easily reduce
to the case $k=1$: assuming that $e,e'$ are incident to $v$ and $S_e
\subset S_{e'}$, we must show that $\{e,e'\}$ is an $S_v$ orbit of the
action of $S_v$ on $\Edges(v)$. If this is not true, then there exists an
edge $e'' \ne e,e'$ incident to $v$ such that $e',e''$ are in the same
$S_v$ orbit (this uses the fact that all orbits of the $S_v$ action on
$\Edges(v)$ have cardinality $\ge 2$). Let $T_{e'}, T_{e''}$ be the
closures of the components of $T-v$ containing $e', e''$, respectively.
Any element of
$S_v$ taking $e'$ to $e''$ restricts to an isomorphism $f \from T_{e'}\to
T_{e''}$. Let $F \from T \to T$ be the isomorphism whose restriction to
$T-(T_{e'} \union T_{e''})$ is the identity, so that $F \restrict
T_{e'}=f$ and $F \restrict T_{e''}=f^\inv$. Then we have $F \in S_e -
S_{e'}$, contradicting that $S_e \subset S_{e'}$, and therefore showing
that $\{e,e'\}$ do form an $S_v$ orbit.
\end{proof}

The condition $S_e = S_{e'}$ is obviously an equivalence relation on
edges, called \emph{stabilizer equivalence}. Condition (c) in the lemma
shows that there are three types of stabilizer equivalence classes. First
is a \emph{singleton}, a class consisting of a single edge $e$; this
occurs when $p(e)$ has no index~2 ends. Second is a \emph{doubleton}, a
pair of edges $e,e'$ sharing an endpoint; this occurs when $p(e)=p(e')$ has
exactly one end of index~2. Third is a \emph{line} in $T$, which occurs
when the image of the line is an edge of $\Gamma$ both of whose ends have
index~2. Condition (a) shows that if $e,e'$ are inequivalent edges then
neither of $S_e, S_{e'}$ is contained in the other.

We shall now define the map $\psi \from T \to T'$. To define $\psi$
on $\Vertices(T)$, note that the map $v \to S_v$ is a bijection between
$\Vertices(T)$ and the maximal compact subgroups of $\Isom T$. Pick a
representative vertex $v$ of each orbit of $\Isom T$; the subgroup
$\Psi(S_v)$ of $\Isom T'$ is compact and must therefore fix some vertex
of $T'$, and we define $\psi(v)$ be any such vertex. Extend
$\psi$ to a $\Psi$-equivariant map $\Vertices(T)\to\Vertices(T')$. Now
extend to a $\Psi$-equivariant map $\psi \from T \to T'$ by mapping each
edge of $T$ to a constant speed geodesic in $T'$.

First we show that $\psi$ is surjective. Since $\psi(T)$ is connected,
nonsurjectivity of $\psi$ would imply that there is a vertex
$v'$ of $T'$ so that some component $C$ of $T'-v'$ is disjoint from
$\psi(T)$. Since $T'$ is thornless, $C$ is unbounded. But this
contradicts cocompactness of $\image(\Psi)$ in $\Isom T'$.

Next we show that $\psi$ is injective on $\Vertices(T)$. If not, consider
$v_1 \ne v_2 \in \Vertices(T)$ such that $\psi(v_1)=\psi(v_2)=w \in
\Vertices(T')$. If $A$ is the closure of the subgroup of $\Isom T$
generated by $S_{v_1} \union S_{v_2}$, then using the fact that $S_{v_1}$
and $S_{v_2}$ are maximal compact subgroups, it follows that $A$ is a
closed, noncompact subgroup of $\Isom T$. Also, $\Psi(A)$ is a closed
subgroup of $\Isom T'$, by properness of $\Psi$. But $\Psi(A)$ is
contained in
$S_w$ which is compact, and so $\Psi(A)$ is compact, contradicting
properness of $\Psi$.

Next we show that for any vertex $v \in \Vertices(T)$, if $v \ne x \in T$
then $\psi(v) \ne \psi(x)$. Arguing by contradiction, suppose
$\psi(v)=\psi(x) = w \in \Vertices(T')$. We may assume $x \in
\interior(e)$ for some edge $e$ of $T$. Arguing as above using properness
of $\Psi$, since $S_w$ is compact it follows that the closure of the
subgroup of $\Isom T$ generated by $S_v \union S_e$ is compact, but
$S_v$ is a maximal compact subgroup and so $S_e \subset S_v$. Let $e_1 *
\cdots * e_k=e$ be the unique embedded edge path in $T$ which starts at
$v$ and ends with $e$. Since $S_e$ stabilizes $v$ it follows that $S_e$
stabilizes each edge in this edge path, and by applying
Lemma~\ref{LemmaIndexTwo} it follows that the edges
$e_1,\ldots,e_k$ are all stabilizer equivalent. Obviously $v$ is not
identified with any point in the interior of $e_1$, and so $k \ge 2$. Let
$L$ be the stabilizer equivalence class of $e_1,\ldots,e_k$, and so either
$k=2$ and $L=e_1*e_2$ is a doubleton, or $L$ is a line in $T$; in either
case we derive a contradiction.

\subparagraph{Case 1: $L = e_1 * e_2$.} Let $S_L$ be the subgroup of
$\Isom T$ stabilizing $L$. The restriction of $S_L$ to $L$ is a standard
$\Z/2$ reflection on an arc, and the map $\psi \from L \to T'$ is $\Z/2$
equivariant. The image $\psi(L)$ is a subtree expressed as a union of two
arcs $\psi(e_1) \union \psi(e_2)$ sharing at least one endpoint
$\psi(v_1)$. By $\Z/2$-equivariance it follows that for $i=1,2$ there
is a subsegment $e'_i$ of $e_i$ incident to $v_i$ such that
$\psi(e'_1)=\psi(e'_2)$, and the sets $\psi(e_1-e'_1)$, $\psi(e_2-e'_2)$
are disjoint from each other and from $\psi(e'_1)=\psi(e'_2)$. Moreover,
$\psi(e_1-e'_1) \ne \emptyset$ if and only if $\psi(e_2-e'_2)
\ne \emptyset$. But this contradicts that $\psi$ identifies $v_0$ with an
interior point of $e_2$.

\subparagraph{Case 2: $L$ is a line.} Say $L = \cdots * e_{-1} * e_0 *
e_1 * \cdots * e_k * e_{k+1} * \cdots$. Let $v_i = e_{i-1} \intersect
e_i$, $i \in \Z$. The restriction to $L$ of the stabilizer of $L$ is a
standard $D_\infinity$ action on the line $L$. The fixed points of the
reflections in $D_\infinity$ are precisely the vertices; let $r_i$ be the
reflection fixing $v_i$. The even vertices $\{v_{2n}\}$ form one orbit
under $D_\infinity$, and the odd vertices $\{v_{2n+1}\}$ form another
orbit. The image $L' = \psi(L)\subset T'$ is the subtree of $T'$ spanned
by $\psi(\Vertices(T))$, the group $D_\infinity$ acts on $L'$, and the
map $L \xrightarrow{\psi} L'$ is $D_\infinity$ equivariant. The action of
$D_\infinity$ on $L'$ is evidently proper and cobounded, and since $L'$
is a tree it follows that there is a $D_\infinity$-invariant line $\ell
\subset L'$ on which the $D_\infinity$ action is standard. Let $w_i \in
\ell$ be the fixed point of $r_i$. Note that for each $w_i$, one of the
following two properties holds, and these properties are equivariant with
respect to $D_\infinity$: either $w_i = \psi(v_i)$; or $\psi(v_i)
\not\in \ell$ and $w_i$ is the closest point to $\psi(v_i)$ on $\ell$. In
either case, it is evident from this description that the arc
$\overline{\psi(v_i) \psi(v_{i+1})}$ contains none of the $w$'s except for
$w_i$ and $w_{i+1}$, and therefore contains none of the $\psi(v)$'s
except for $\psi(v_i)$ and $\psi(v_{i+1})$. But this contradicts that
$\psi(v_0)$ lies on $\overline{\psi(v_{k-1}) \psi(v_k)}$.

The argument in the last paragraph gives a little more information: it
shows that for any stabilizer equivalence class which is a line $L$, any
two nonadjacent edges of $L$ have disjoint images under $\psi$.

\bigskip

We have shown that $\psi \from T \to T'$ does not identify any vertex of
$T$ with any other point of $T$.

Next we show, for any two nonadjacent edges $e,e'$ of $T$, that
$f(\interior(e)) \intersect f(\interior(e')) = \emptyset$. Suppose not:
there exists $x \in \interior(e)$, $x' \in \interior(e')$ such that
$f(x)=f(x')=y$. Letting $A$ be the closure of the subgroup of $\Isom T$
generated by $S_e \union S_{e'}$, it follows that $\Psi(A)$ stabilizes the
point $y$. Using properness of $\Psi$ it follows that $A$ is compact.
This implies that $A$ stabilizes some vertex $v$ of $T$, and so $S_e
\union S_{e'} \subset A \subset S_v$. Let $e=e_0 * \cdots * e_k$ be the
shortest edge path from $e$ to $v$, and similarly for $e'=e'_0 * \cdots *
e'_{k'}$. It follows that $S_e \subset S_{e_i}$, $i=0,\ldots,k$ and
$S_{e'} \subset S_{e'_{i'}}$, $i' = 0,\ldots,k'$. Applying
Lemma~\ref{LemmaIndexTwo} it follows that the stabilizer equivalence
classes $L$ and $L'$ of $e$ and $e'$ contain $e_0,\ldots,e_k$ and
$e'_0,\ldots,e'_{k'}$, respectively, and so $L$ and $L'$ are both
incident to the vertex $v$. Also, we already know that
$e,e'$ are not stabilizer equivalent to each other because disjoint edges
of a stabilizer equivalence class have disjoint images, and so $L
\ne L'$. Using the description above of a stabilizer equivalence class and
its image under $\psi$, and using the fact that $\interior(e)$ and
$\interior(e')$ are disjoint from $\psi(\Vertices(T))$, we may reduce to
the case that $k,k' \le 1$; and since $e,e'$ are not adjacent at least one
of $k,k'$ is $=1$. Consider the case where one of $k,k'$ equals $0$, say
$k=1, k'=0$ (the other case, where $k'=k=1$, is similar and is left to the
reader). Then there must be $x_1 \in \interior(e_1)$ such that
$\psi(x)=\psi(x_1)=\psi(x') = y$ in $T'$. Let $w = e_0 \intersect e_1$ and
choose $g \in S_w$ which interchanges $e=e_0$ with $e_1$, and so $w$
interchanges
$x$ with $x_1$. By equivariance under $g$, the edge $e''=g(e')$
contains a point $x''$ such that $\psi(x'')=y$, and by the argument just
given it follows that the stabilizer equivalence classes $L', L''$ of
$e', e''$ are distinct but are adjacent to a common vertex. But this
is impossible, because $e', e''$ are separated from each other by the
edges $e_0,e_1$ of the line $L$, and the lines $L, L', L''$ are distinct
stabilizer equivalence classes in $T$.

Next we show that each edge $e$ of $T$ contains a point denoted $m_e$
such that $\psi(x) \ne \psi(m_e)$ for any $x \ne m_e$; we may choose the
points $m_e$ equivariantly with respect to $\Isom T$. The point $m_e$ is
called the \emph{midpoint} of $e$ and the closures of the
two components of $e-m_e$ are called the \emph{halves} of $e$. To see why
$m_e$ exists, let $v,w$ be the endpoints of $e$. Let $e_v$ be the longest
subsegment of $e$ which is identified via $\psi$ with a subsegment of
another edge incident to $v$, and similarly for $e_w$. It follows that
$e_v \intersect e_w = \emptyset$ for otherwise there would be edges
$e',e'' \ne e$ incident to $v,w$ such that $\psi(e') \intersect \psi(e'')
\ne \emptyset$, contradiction. We can then take $m_e$ to be any point of
$\interior(e) - (e_v \union e_w)$.

Now we may give a global description of the map $\psi \from T \to T'$.
Given a vertex $v$, define $\Star(v)$ to be the union, over all edges
$e$ incident to $v$, of the half of $e$ containing $v$. Note that
$\Star(v)$ is invariant under $S_v$. The restriction of $\psi$ to
$\Star(v)$ is an $S_v$-equivariant family of partial Stallings folds
\cite{Stallings:folds}: the half-edges forming $\Star(v)$ are subdivided
and then folded, no two half-edges being entirely folded together. These
are the only identifications made by $\psi$. Note that for any half-edge
$e$ incident to $v$, if $e$ is partially folded with any other half-edge
then all half-edges in the $S_v$ orbit of $e$ are partially folded
together to form a single path in $T'$; these paths, one for each
partially folded orbit of $v$-half-edges, may then undergo further partial
foldings among each other.

The description of the map $\psi \from T \to T'$ shows that the
edge-indexed graph $\hat\Gamma = T' / \Psi(\Isom T)$ is a blowup of
$\Gamma$: for any $v \in \Vertices(T)$, the partial folds performed on
$\Star(v)$ are represented downstairs in $\Gamma$ by a blowup of the
vertex $p(v)$ of $\Gamma$, and doing this for each vertex of $\Gamma$ we
obtain a blowup $\Gamma\leftarrow\hat\Gamma$. The induced map $\hat\Gamma =
T' / \Psi(\Isom T) \to T' / \Isom T' = \Gamma'$ is a covering map, by
Lemma~\ref{LemmaCovering}. By hypothesis, $\Gamma$ has no blowup and
proper subcover, and so $\hat\Gamma$ is isomorphic to $\Gamma'$. Also by
hypothesis, $\Gamma'$ has no index~1 edges, and so the blowup is trivial
and the map $\psi \from T \to T'$ is an isometry.

This completes the proof of Proposition~\ref{PropositionIsometry}. 
\end{proof} 

Combining Proposition~\ref{PropositionIsometry}
with Theorem~\ref{TheoremMaximal} we immediately have:

\begin{corollary}
\label{CorollaryMaximal}
A bounded valence, bushy, index~1 normalized tree $T$ is
maximally symmetric if and only if its quotient edge-indexed graph
$\Gamma=T/\Isom T$ has no blowup and proper subcover. It follows that
the correspondence taking a tree $T$ to its quotient graph $T/\Isom T$
sets up a bijection between isometry classes of index~1 normalized
maximally symmetric trees and isomorphism classes of finite edge-indexed
graphs with no index~1 edge and with no blowup and proper subcover. More
generally, a bounded valence, bushy, cocompact, thornless tree $T$ is
maximally symmetric if and only if the edge-indexed graph
$\Gamma=T/\Isom T$ satisfies the following: each index~1 edge is of
index~1--1; the collection of these edges forms a forest $F$; and the
collapsed graph $\Gamma/F$ has no blowup and proper subcover.
\qed\end{corollary}

\subsection{Proof of Theorem~\ref{TheoremModelGeometries} and
Corollary~\ref{CorollaryVirtuallyFree}}

Fix a bounded
valence, bushy tree $\tau$. Suppose that $\G$ is a uniform, cobounded
subgroup of $\tau$. The inclusion map $\G \to
\QI(\tau)$ factors through a map $\alpha \from \G \to \QIhat(\tau)$ which
is a cobounded quasi-action of $\G$ on $\tau$. Applying
Theorem~\ref{TheoremTreeQuasiActions} there is a bounded valence, bushy,
cocompact tree $T$, a cobounded action $\phi \from
\G \to \Isom T$, and a quasi-conjugacy $f \from \tau \to T$ from the
quasi-action $\alpha$ to the action $\G$. By trimming thorns we may assume
that $T$ is thornless, and so the natural homomorphism $\Isom T\to
\QI(T)$ is an embedding. Consider the edge-indexed graph $\Gamma = T /
\Isom T$. Applying Corollary~\ref{CorollaryPumpingUpAlgorithm} there is a
collapse and subcover $\Gamma \xrightarrow{p} \Gamma_1$, and a blowup and
subcover
$\Gamma_1 \xleftarrow{q} \Gamma_2 \xrightarrow{\mu}
\Gamma_3$, such that $\Gamma_3$ has no index~1 edge and no blowup and
proper subcover. The maps $p,q,\mu$ lift to quasi-isometries of universal
covering trees $T \xrightarrow{P} T_1$, $T_1 \xleftarrow{Q} T_2$, $T_2
\xrightarrow{M} T_3$. Applying Corollary~\ref{CorollaryMaximal} the tree
$T_3$ is maximally symmetric, and applying Theorem~\ref{TheoremMaximal}
the group $\Isom T_3$ is a maximal uniform subgroup of $\QI(T_3)$. The
graphs $\Gamma$, $\Gamma_1$, and $\Gamma_3$ are their own minimal
subcovers, and so $D(\Gamma)$, $D(\Gamma_1)$, and $D(\Gamma_3)$ equal
$\Isom T$, $\Isom T_1$, and $\Isom T_3$ respectively. We have
$\ad_f(\G) \subgroup \Isom T$, $\ad_P(\Isom T) \subgroup \Isom T_1$,
$\ad_Q^\inv(\Isom T_1)=D(\Gamma_2)$, and $\ad_M(D(\Gamma_2)) \subgroup
\Isom T_3$. Letting $F = M \composed \bar Q \composed P \composed f
\from \tau \to T_3$, where $\bar Q \from T_1 \to T_2$ is a coarse inverse
of $Q$, it follows that $\ad_F = \ad_M \composed \ad_Q^\inv \composed \ad_P
\composed\ad_f\from \QI(\tau) \to \QI(T_3)$ is an isomorphism inducing a
bijection between uniform subgroups, and $\ad_F^\inv(\Isom T_3)$ is a
maximal uniform subgroup of $\QI(\tau)$ containing~$\G$. This proves (1).
Moreover, if $\G$ is a maximal uniform cobounded subgroup of $\tau$ it
follows that $\G = \ad_F^\inv(\Isom T_3)$, proving the first sentence of
(2). 

To prove the second sentence of (2), suppose that $\G =
\ad_{f'}(\Isom T')$ where $T'$ is a bounded valence, bushy, cocompact
tree and $\Gamma' = T' / \Isom T'$ has no index~1 edges, and where $f'
\from T' \to\tau$ is a quasi-isometry. We thus have a quasi-isometry $T_3
\xrightarrow{H=f'{}^\inv \composed F} T'$ with the property that
$\ad_H$ takes $\Isom T_3$ into $\Isom T'$ by a continuous, proper,
cocompact monomorphism. Applying Proposition~\ref{PropositionIsometry} it
follows that $\ad_H = \ad_{H'}$ for some isometry $H' \from T_3 \to T'$,
proving (2).

To prove part (3) of the theorem, suppose that $\G,\G'$ are maximal
uniform subgroups of $\QI(\tau)$ and that $F \G F^\inv = \G'$ for some $F
\in \QI(\tau)$. Applying part (2) we have $\G=\ad_f(\Isom T)$, $\G' =
\ad_{f'}(\Isom T')$ for some maximally symmetric trees $T$ and
quasi-isometries $f\from T\to\tau$, $f'\from T' \to \tau$. Also, part (2)
shows that $T,T'$ are uniquely determined up to isometry by $\G,\G'$.
Moreover, part (2) shows that $T,T'$ are isometric to each other, because
the isomorphism $\ad_{f'}^\inv \composed ad_F \composed \ad_f \from
\Isom(T)\to \Isom(T')$ is equal to $\ad_h$ for some isometry $h\from T \to
T'$. Thus, to each conjugacy class of maximal uniform subgroups of
$\QI(\tau)$ there corresponds a well-defined isometry class of maximally
symmetric trees. This correspondence is a surjection, because for every
maximally symmetric tree $T$ there exists a quasi-isometry $f \from T \to
\tau$, and so $\G = \ad_f(\Isom T)$ is a maximal uniform subgroup of
$\QI(\tau)$ by Theorem~\ref{TheoremMaximal}. Also this correspondence is an
injection, for suppose we have maximal uniform subgroups $\G=\ad_f(\Isom
T), \G'=\ad_{f'}(\Isom T')$  where $T,T'$ are maximally symmetric, $f \from
T \to \tau, f' \from T' \to\tau$ are quasi-isometries, and $T,T'$ are
isometric; choosing an isometry $h \from T \to T'$ it follows that $F=[f'
\composed h \composed f] \in \QI(\tau)$ conjugates $\G$ to $\G'$.

We can prove part (4) of Theorem~\ref{TheoremModelGeometries} by simply
noting some examples of edge-indexed graphs with no index~1 edge nor any
blowup and proper subcover. 

For unimodular examples, the edge-indexed graphs 
$
\xymatrix{
\bullet \ar@{-}[r]^<<{p}^>>{q} & \bullet 
}
$
with $p > q \ge 2$ clearly cannot be pumped up: they have no proper
subcovers; and they have no blowups. Their universal
covers give countably many isometry classes of maximally symmetric
unimodular trees.

For nonunimodular examples, consider edge-indexed graphs of the form
$$
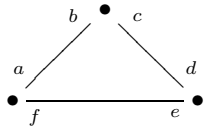
\xymatrix{
& \bullet \ar@{-}[dl]_>>{a}_<<{b} \ar@{-}[dr]^<<{c}^>>{d} \\
\bullet \ar@{-}[rr]_>>{e}_<<{f} && \bullet
}
$$
Assign integer values $\ge 2$ to $a,b,c,d,e,f$ so that the numbers $a+f,
b+c, d+e, a+1, b+1, c+1, d+1, e+1, f+1$ are pairwise unequal. These
numbers are the total indices that can occur for vertices in any
blowup, and it follows that each blown up graph has vertices with distinct
total indices and therefore has no proper subcover.
Clearly we may make the choices so that $\frac{bdf}{ace} \ne 1$ and so the
edge-indexed graph is not unimodular.

See Proposition~\ref{PropositionInfinite} below for a more satisfactory
enumeration of examples.

This completes the proof of Theorem~\ref{TheoremModelGeometries}. \qed

\begin{proof}[Proof of Corollary~\ref{CorollaryVirtuallyFree}] Let $G$ be
a virtually free group of finite rank~$\ge 2$. There is a finite graph of
finite groups $\Gamma'$ whose fundamental group is isomorphic to $G$
\cite{KPS:virtuallyfree}, and so $G$ acts properly discontinuously and
cocompactly on the Bass-Serre tree $\tau$ of $\Gamma'$. The kernel of the
action $\phi \from G \to \Isom\tau$ is finite.
By Theorem~\ref{TheoremModelGeometries}, the group $\phi(G)$ is contained
in a maximal uniform subgroup of $\QI(\tau)$, and that subgroup is equal to
$\ad_f(\Isom T)$ for some maximally symmetric tree $T$ and some
quasi-isometry $f \from T \to \tau$. The action $\ad_f^\inv \composed \phi
\from G \to \Isom T$ is properly discontinuous and cocompact, and the
graph of groups $\Gamma = T / \ad_f^\inv\composed\phi(G)$ has fundamental
group isomorphic to $G$ and Bass-Serre tree $T$.
\end{proof}

Finally, we improve upon
Theorem~\ref{TheoremModelGeometries} part (4) as follows:

\begin{proposition}
\label{PropositionInfinite}
Let $\Gamma$ be a finite, connected graph with no loops and no bigons and
with at least one edge. Then $\Gamma$ has infinitely many distinct
edge-indexings with no index~1 edge and no blowup and
proper subcover. Moreover, infinitely many of them are unimodular, and if
$\Gamma$ is not a finite tree then infinitely many of them are
nonunimodular.
\end{proposition}

\begin{proof} One way to proceed with the proof is to construct
sufficiently many examples, that is, to describe some special scheme for
constructing sufficiently many edge-indexings of $\Gamma$ with the desired
properties. Instead we shall give a general scheme for enumerating all of
the appropriate edge indexings of $\Gamma$, in effect enumerating
the maximally symmetric trees $T$ with $T / \Isom(T)$ isomorphic to
$\Gamma$. From the description, the infinitude of appropriate
edge-indexings of $\Gamma$ will follow. 

The enumeration scheme is carried out completely for one example $\Gamma$,
after the conclusion of the proof. The reader may want to refer to this
example while perusing the proof.

Instead of indexing the ends of $\Gamma$ with actual
numerical values, index them with variables $x_1,\ldots,x_N$. An actual
edge-indexing of $\Gamma$ without index~1 ends corresponds to an assignment
of $(x_1,\ldots,x_N) \in \{2,3,\ldots\}^N$. 

We shall construct $X \subset \R^N$, a finite union of affine subspaces of
$\R^N$ defined over $\rational$, such that an edge-indexing
$(x_1,\ldots,x_N)$ is in $X$ if and only if there exists a blowup and
proper subcover for the edge indexing $(x_1,\ldots,x_N)$. Then, when
$\Gamma$ is not a finite tree, we shall construct a certain homogeneous
subvariety $Y\subset \R^N$ of degree~$\ge 3$ defined over $\rational$, not
contained in any degree~1 subvariety except for all of $\R^N$, such that an
edge-indexing $(x_1,\ldots,x_N)$ is in $Y$ if and only if
$(x_1,\ldots,x_N)$ is a unimodular. The conclusions of the theorem will
quickly follow from the form of $Y$.

Let $\Gamma \xleftarrow{F_1} \Gamma_1,\ldots,
\Gamma\xleftarrow{F_K}\Gamma_K$ denote all possible blowups of $\Gamma$,
where each end of $\Gamma_k$ is indexed with one of the variables
$x_1,\ldots,x_N$ or with the integer $1$, as follows. As usual $F_k$
denotes a subtree of $\Gamma_k$ consisting of edges of index 1--1, the
cellular map $\Gamma_k \to \Gamma$ collapses each component of $F_k$ to a
point and is otherwise one-to-one, and under this collapse there is a
bijection between $\Ends(\Gamma_k)-\Ends(F_k)$ and $\Ends(\Gamma)$. This
bijection is used to index $\Ends(\Gamma_k)-\Ends(F_k)$ with the variable
edge indices $x_1,\ldots,x_N$. We assume that one of these blowups, say
$\Gamma\leftarrow\Gamma_1$, is actually the identity map on $\Gamma$, and
so $F_1=\emptyset$. 

Given $k=1,\ldots,K$, let $\mu_{kj} \from\Gamma_k \to
\Gamma_{kj}$, $j=1,\ldots,J(k)$ denote all the proper subcovers of
$\Gamma_k$, where $\Gamma_{kj}$ is a graph with a variable $y_e$ indexing
each end $e \in \Ends(\Gamma_{kj})$. The even covering equations for
$\mu_{kj}$ form a system $\E_{kj}$ of first degree equations in the
variables $\{x_1,\ldots,x_N\}\union\{y_e\suchthat
e\in\Ends(\Gamma_{kj})\}$. 

We examine the even covering system $\E_{kj}$ more carefully. To each end
$e\in\Ends(\Gamma_{kj})$ located at a vertex $v \in
\Vertices(\Gamma_{kj})$, and to each vertex $w \in
\Vertices(\Gamma_k) \intersect\mu_{kj}^\inv(v)$, there  corresponds an
even covering equation whose right hand side is
$y_e$ and whose left hand side is a sum of those variables
$x_1,\ldots,x_N$ labelling ends of $\mu_{kj}^\inv(e) \intersect \Ends(w)$
plus an integer equal to the cardinality of
$\Ends(F_k) \intersect \mu_{kj}^\inv(e) \intersect \Ends(w)$; let
$\E_{kj}^{ew}$ denote this equation. Let $\E_{kj}^{e}$ be the system of
equations $\E_{kj}^{ew}$, $w\in \mu_{kj}^\inv(v)$; and let $\E_{kj}$ be
the system of equations $\E_{kj}^e$, $e \in \Ends(\Gamma_{kj})$. 

If the system $\E_{kj}$ is inconsistent then we may discard $\Gamma
\leftarrow \Gamma_k \rightarrow \Gamma_{kj}$ as a candidate blowup and
proper subcover. We may characterize inconsistency of $\E_{kj}$ as
follows. Note that for $e\ne e'$ the subsystems $\E_{kj}^e$ and
$\E_{kj}^{e'}$ have disjoint variable sets, and so the system $\E_{kj}$ is
consistent if and only if each of the subsystems
$\E_{kj}^e$ is consistent. The subsystem $\E_{kj}^e$ is inconsistent if and
only if there exist $w,w' \in\mu_{kj}^\inv(v)$ such that each of the two
sets $\mu_{kj}^\inv(e)\intersect\Ends(w)$,
$\mu_{kj}^\inv(e)\intersect\Ends(w')$ lies entirely in
$\Ends(F_k)$ but these two sets have different cardinalities. 

Assume now that the system $\E_{kj}$ is consistent. In each subsystem
$\E^e_{kj}$, since the variable $y_e$ occurs alone on the right hand side
of each equation in $\E^e_{kj}$, we may eliminate $y_e$; after this
elimination, if there are any equations of the form
(constant)$=$(constant) we may eliminate them as well. This produces a new
system of equations $\bar\E^e_{kj}$ in the variables $x_1,\ldots,x_N$. The
system $\bar\E_{kj}$ is the union of the systems
$\bar\E^e_{kj}$ for $e \in \Ends(\Gamma_{kj})$; we call $\bar\E_{kj}$ the
system of \emph{reduced even covering equations} for $\mu_{kj}$. Let
$X_{kj} \subset \R^N$ denote the solution set of $\bar\E_{kj}$, an affine
subspace of $\R^N$ defined over $\rational$. 

We claim that if $\bar\E_{kj}$ is vacuous, meaning that $X_{kj} =
\reals^N$, then $\mu_{kj} \from \Gamma_k \to \Gamma_{kj}$ is a graph
isomorphism and so is \emph{not} a proper subcover. 

To prove the claim, assume that $\bar\E_{kj}$ is vacuous. Since the
variables in distinct subsystems $\bar\E^e_{kj}$ are distinct, it follows
that each $\bar\E^e_{kj}$ is vacuous. Let $v$ be the vertex incident to
$e$. Vacuity of $\bar\E^e_{kj}$ implies one of two possibilities. In
the first possibility, $\mu_{kj}^\inv(v) \subset F_k$; in this case each
equation $\E^{ew}_{kj}$ has the form (constant)$= y_e$, and so elimination
removes all these equations. In the second possibility, $\mu_{kj}^\inv(v)$
is a single vertex $w$; in this case $\E^e_{kj}$ consists of the single
equation $\E^{ew}_{kj}$, with right hand side $y_e$, and elimination of
$y_e$ removes this equation. In all other cases,
$\bar\E^e_{kj}$ is nonvacuous.

Partition the set $\Vertices(\Gamma_k)$ into two sets: $V$ consists of all
vertices incident to some variably indexed end; and $V'$ consist of all
the rest, namely those vertices incident only to ends of index~1, the
interior vertices of $F_k$. It follows from the previous paragraph that
$\mu_{kj}(V)\intersect\mu_{kj}(V') = \emptyset$ and that $\mu_{kj}$ is
one-to-one on $V$. Since $\Gamma_k$ has neither loops nor bigons it
follows further that $\mu_{kj}$ is one-to-one on the union of all edges
with both ends in the set $V$. Consider now the restriction of $\mu_{kj}$
to $F_k$; we have already seen that $\mu_{kj}$ is one-to-one on $F_k
\intersect V$ and the remaining vertices of $F_k$ are mapped disjointly
from $V$. But this implies that $\mu_{kj}$ is one-to-one on vertices and
edges of $F_k$. Thus we see that $\mu_{kj}$ is an isomorphism, proving the
claim.

To summarize, we have shown that the edge-indexings of $\Gamma$ which
\emph{do} have a blowup and proper subcover are those which lie on a finite
union $X=\Union_{kj} X_{kj}$ of codimension~$\ge 1$ affine subspaces of
$\reals^N$ defined over $\rational$. It follows that there are infinitely
many edge indexings with integer values~$\ge 2$ that do not lie in $X$,
and so infinitely many edge-indexings with no index~1 edge and no blowup
and proper subcover.  If
$\Gamma$ is a finite tree then this finishes the proof, because every
edge-indexing of $\Gamma$ is unimodular.

Suppose now that $\Gamma$ is not a finite tree, but instead is a graph of
rank $R$. Let $c_1,\ldots,c_R$ be simple, closed, oriented edge paths
whose corresponding 1-cycles give a basis for $H_1(\Gamma)$. Let $\ell(r)
\ge 3$ be the number of edges in $c_r$. Applying the canonical cocycle
$\xi$ to $c_r$ we obtain an equation of the form
$$\xi(c_r) = \frac{x_{n_1} \cdot \ldots \cdot
x_{n_{\ell(r)}}}{x_{n_{\ell(r)+1}} \cdot \ldots \cdot x_{n_{2\ell(r)}}}
$$
where all $2\ell(r)$ of the variables are distinct.
Setting this equal to $1$ and clearing the denominator we thus obtain the
following homogeneous equation of degree $\ell(r)$:
$$x_{n_1} \cdot \ldots \cdot
x_{n_{\ell(r)}} = x_{n_{\ell(r)+1}} \cdot \ldots \cdot x_{n_{2\ell(r)}}
$$
Let $Y$ be the simultaneous solution variety of this system of equations
for $r=1,\ldots,R$; the points on $Y$ with integer coordinates $\ge 1$ are
precisely the unimodular edge-indexings of $\Gamma$. Note that rational
points are dense in $Y$.

Clearly $Y$ has codimension~$\ge 1$ and so there are infinitely many
edge-indexings in the complement of $X \union Y$, that is, infinitely many
nonunimodular edge-indexings which have no index~1 edge and no blowup and
proper subcover.

For the unimodular case, note that the homogeneous variety $Y$ is not of
degree~1, indeed $Y$ is not contained in any linear subspace of
$\R^N$ except for $\R^N$ itself. To see why we rewrite the
defining equations for $Y$ as follows. Choosing a maximal tree $T$ of
$\Gamma$ we may push all variables in $T$ to the right hand side and all
variables not in $T$ to the left hand side, and then reorder the
variables, obtaining a set of defining equations for
$Y$ of the form
\begin{align*}
\frac{x_1}{x_2} &= f_1(x_{2R+1},\ldots,x_N) \\
  & \vdots \\
\frac{x_{2R-1}}{x_{2R}} &= f_R(x_{2R+1},\ldots,x_N)
\end{align*}
where each $f_r$ is a quotient of homogeneous monomials of equal degree
$\ell(r)-1 \ge 2$ with no variable occurring more than once in $f_r$. From
this it is obvious that $Y$ is not contained in any proper linear subspace
of $\R^N$. 

Decompose $X$ as $X=X' \union X''$ where $X'$ is the union of those
$X_{kj}$ which are homogeneous and $X''$ is the union of those which are
not homogeneous. As we have just seen, $Y \not\subset X'$, and it follows
that the rational rays lying in $Y-X'$ are dense in
$Y$. For each such ray $\rho$ the intersection $\rho \intersect X''$ is
finite, and so assuming that $\rho$ points into the positive orthant of
$\R^N$ it follows that the number of integer points on $\rho - X''$ with
coordinates $\ge 2$ is infinite.
\end{proof}

\break
\paragraph{An example}
We close with a consideration of the edge-indexed graph 
$$
\Gamma \quad = \quad \xymatrix{
\bullet \ar@{-}[rr]^<<<<{a}^>>{b} && \bullet \ar@{-}[rr]^<<{c}^>>>>{d}
&&
\bullet }
$$
We enumerate the thirteen different blowup and proper subcovers of
$\Gamma$, writing down their reduced even covering equations. We thus
obtain $X$ as a union of thirteen affine subspaces of $(a,b,c,d)$-space,
although it turns out that six subspaces suffice.

To help in the enumeration, consider an edge-indexed graph which is an
arc with $m$ edges. A subcover \emph{without subdivision} of such a graph
must also be an arc, because vertex valence cannot
increase under a subcover. Moreover, a covering map (without subdivision)
from an $m$-edged arc to an $n$-edged arc must be an $m/n$-folding map,
under which the $m$ edges of the domain are partitioned into $m/n$ subarcs
each with $n$-edges, each mapped by a graph isomorphism to the range; this
follows because the ends located at a vertex of the domain must map
surjectively to the ends located at the image vertex of the range. A
subcover with subdivision is similarly described, with the proviso that
each subdivision point must be a fold point.

We start by enumerating the proper subcovers of $\Gamma$ itself, by
choosing a subdivision followed by a fold.
\begin{itemize}
\item[(1)] Bifold $\Gamma$ over a single edge; the reduced even covering
equation is $a=d$. 
\item[(2)] Subdivide the $a$--$b$ edge and then trifold over a single
edge; the equations are \break $a=b+c, d=2$.
\item[(3)] Subdivide the $c$--$d$ edge and trifold; the equations are
$a=2, d=b+c$. 
\item[(4)] Subdivide both edges and quadrifold; the equations are
$a=b+c=d$.
\end{itemize}
Next there are the proper subcovers of the unique nontrivial blowup
$\Gamma'$ of $\Gamma$:
$$
\Gamma' \quad=\quad \xymatrix{
\bullet \ar@{-}[r]^<<{a}^>>{b} & 
   \bullet \ar@{-}[r]^<<{1}^>>{1} &
   \bullet \ar@{-}[r]^<<{c}^>>{d} & 
   \bullet
}
$$
The proper subcovers of $\Gamma'$ are:
\begin{itemize}
\item[(5)] Trifold of $\Gamma'$; the equations are $a = c+1, d
= b+1$. This subcover was depicted earlier for $(a,b,c,d)=(4,5,3,6)$. 
\item[(6)] Subdivide the $a$--$b$ edge and quadrifold; the equations are
$a=b+1=d, c=1$.
\item[(7)] Subdivide the 1--1 edge and quadrifold; the equations are
$a=d=2, b=c$.
\item[(8)] Subdivide the 1--1 edge and bifold; the equations are $a=d,
b=c$.
\item[(9)] Subdivide the $c$--$d$ edge and quadrifold; the equations are
$a=c+1=d, b=1$.
\item[(10)] Subdivide the $a$--$b$ and 1--1 edges and pentafold; the
equations are $a=b+1=c+1, d=2$.
\item[(11)] Subdivide the $a$--$b$ and $c$--$d$ edges and pentafold; the
equations are $a=b+1=2, d=c+1=2$, or equivalently, $a=d=2, b=c=1$.
\item[(12)] Subdivide the 1--1 and $c$--$d$ edges and pentafold; the
equations are $a=2, d=c+1=b+1$.
\item[(13)] Subdivide all three edges and hexafold; the equations are
$a=b+1=c+1=d$.
\end{itemize}
Let $X_{(i)}$ be the affine subspace defined by case $(i)$. Then clearly
$X_{(1)} = \{a=d\}$ contains $X_{(i)}$ for $i=4,6,7,8,9,11,13$.
Thus we may eliminate the latter seven equations, and we obtain the
following set of six affine subspaces whose union is $X$:
\begin{align*}
X_{(1)} &= \{a = d\} \\
X_{(2)} &= \{a = b+c,  \quad d=2\} \\
X_{(3)} &= \{a = 2,  \quad d=b+c\} \\
X_{(5)} &= \{a = c+1,  \quad d=b+1\} \\
X_{(10)} &= \{a = b+1 = c+1,  \quad d=2\} \\
X_{(12)} &= \{a = 2, \quad d=c+1=b+1\}
\end{align*}
The edge-indexings of $\Gamma$ which correspond to maximally symmetric
trees are precisely the quadruples $(a,b,c,d)$ of integers $\ge 2$ which
lie on none of these six subspaces.


\providecommand{\bysame}{\leavevmode\hbox to3em{\hrulefill}\thinspace}

\bigskip\noindent
\textsc{\scriptsize
Lee Mosher:\\
Department of Mathematics,
Rutgers University,
Newark, NJ 07102\\
mosher@andromeda.rutgers.edu
}

\bigskip\noindent
\textsc{\scriptsize
Michah Sageev:\\
Technion, Israel University of Technology,
Dept.\ of Mathematics,
Haifa 32000, Israel\\
sageevm@techunix.technion.ac.il
}

\bigskip\noindent
\textsc{\scriptsize
Kevin Whyte:\\
Dept.\ of Mathematics,
University of Chicago,
5734 S.\ University Ave.,
Chicago, IL 60637\\
kwhyte@math.uchicago.edu
}

\vfill

\end{document}